\documentclass[sn-mathphys-num]{sn-jnl}

\usepackage{graphicx}%
\usepackage{multirow}%
\usepackage{amsmath,amssymb,amsfonts}%
\usepackage{amsthm}%
\usepackage{mathrsfs}%
\usepackage[title]{appendix}%
\usepackage{xcolor}%
\usepackage{textcomp}%
\usepackage{manyfoot}%
\usepackage{booktabs}%
\usepackage{algorithm}%
\usepackage{algorithmic}
\usepackage{listings}%
\usepackage{cleveref}
\usepackage{tikz} 
\usepackage{marginnote}

\theoremstyle{thmstyleone}%

\theoremstyle{thmstyletwo}%

\theoremstyle{thmstylethree}%

\begin{document}

\title[Randomized and Inner Product-free Krylov Methods]{Randomized and Inner-product Free Krylov Methods for Large-scale Inverse Problems}

\author[1]{\fnm{Malena} \sur{Sabat\'{e} Landman}}\email{Malena.SabateLandman@maths.ox.ac.uk}

\author[2]{\fnm{Ariana N.} \sur{Brown}}\email{ariana.brown@emory.edu}

\author[3]{\fnm{Julianne} \sur{Chung}}\email{jmchung@emory.edu}

\author*[4]{\fnm{James G.} \sur{Nagy}}\email{jnagy@emory.edu}

\affil[1]{\orgdiv{Mathematical Institute}, \orgname{University of Oxford}, \orgaddress{\city{Oxford}}, \postcode{OX2 6GG}, \country{UK}}

\affil[2]{\orgdiv{Dept.~of Mathematics}, \orgname{Emory University}, \orgaddress{\city{Atlanta}}, \postcode{30322}, \state{GA}, \country{USA}}

\affil[3]{\orgdiv{Dep.~of Mathematics}, \orgname{Emory University}, \orgaddress{\city{Atlanta}}, \postcode{30322}, \state{GA}, \country{USA}}

\affil[4]{\orgdiv{Dep.~of Mathematics}, \orgname{Emory University}, \orgaddress{\city{Atlanta}}, \postcode{30322}, \state{GA}, \country{USA}}

\abstract{
Iterative Krylov projection methods have become widely used for solving large-scale linear inverse problems. However, methods based on orthogonality include the computation of inner-products, which become costly when the number of iterations is high; are a bottleneck for parallelization; and can cause the algorithms to break down in low precision due to information loss in the  projections. 
Recent works on inner-product free Krylov iterative algorithms alleviate these concerns, but they are quasi-minimal residual rather than minimal residual methods. This is a potential concern for inverse problems where the residual norm provides critical information from the observations via the likelihood function, and we do not have any way of controlling how close the quasi-norm is from the norm we want to minimize. In this work, we introduce a new Krylov method that is both inner-product-free and minimizes a functional that is theoretically closer to the residual norm.  
The proposed scheme combines an inner-product free Hessenberg projection approach for generating a solution subspace with a randomized sketch-and-solve approach for solving the resulting strongly overdetermined projected least-squares problem. Numerical results show that the proposed algorithm can solve large-scale inverse problems efficiently and without requiring inner-products.
}


\maketitle

\section{Introduction}\label{sec1}

Inverse problems arise in many different applications including medical and geophysical imaging, electromagnetic scattering, machine learning, and image deblurring \cite{Chung2011Inverse,Hansen2010,Vogel2002,Zhdanov2002}. 
Typically, the goal of an inverse problem is to estimate some unknown quantities or parameters of a system, given observations or indirect measurements (e.g., taken on the exterior of the object).
In particular, we consider linear discrete inverse problems of the form,
\begin{equation} \label{eq:ip} 
b = Ax_{\rm true} + e, 
\end{equation}
where $A \in \mathbb{R}^{m \times n}$ models the forward problem, $x_{\rm true} \in \mathbb{R}^{n}$ is the unknown solution we want to approximate, $b \in \mathbb{R}^{m}$ is the vector of observed data, and $e \in \mathbb{R}^{m}$ represents noise and other measurement errors. Given $b$ and $A$, the goal is to estimate $x_{\rm true}$, but there are various computational challenges.

For many applications, the number of unknowns $n$ may be very large and the forward model matrix $A$ (and its adjoint) can only be accessed via matrix-vector multiplications.  Thus, iterative methods are often used to compute approximations of $x_{\rm true}$.  Moreover, these inverse problems are usually ill-posed in the sense that small perturbations in the observation can cause large perturbations in the computed solution. This is due to the fact that the singular values of $A$ decay and cluster at zero without a gap between consecutive values, and 
because the singular vectors corresponding to small singular values tend to be highly oscillatory. 
Indeed, the noiseless right-hand side $Ax_{\rm true}$ satisfies the discrete Picard condition \cite{Hansen2010}, while the noise is captured mostly in the singular subspace associated to the small singular values. Moreover, due to the presence of noise in the measurements, $b$ might not be in the range of $A$. Regularization can be used to address ill-posedness and can take many forms. In iterative regularization, standard iterative solvers (e.g., Krylov methods) are applied to the least squares problem,
\begin{equation} \label{eq:LS} 
\min_{x \in \mathbb{R}^{n}} \|Ax-b\|_2, 
\end{equation}
and regularization is obtained by early stopping.  Variational regularization is another form of regularization, where prior knowledge about the solution properties (e.g., smoothness or sparsity) is included in the regularization term.  We focus on standard Tikhonov regularization,
\begin{equation} \label{eq:LS2} 
\min_{x \in \mathbb{R}^{n}} \|Ax-b\|_2^2 +\lambda \|x\|_2^2,
\end{equation}
where $\lambda>0$ is a regularization parameter that balances the effect of the regularization against the fitting of the data to the noisy measurements (i.e., the likelihood function).

There exists a myriad of iterative methods to solve minimization problems \eqref{eq:LS} and \eqref{eq:LS2} efficiently. However, the capacity we have to measure and store data is constantly increasing, pushing the limits of traditional least-squares solvers in terms of speed and memory requirements. Recently, different directions have emerged to tackle more challenging scenarios.  Some (partial) solutions are motivated in large part by the evolution of commercial hardware.  For example, new methods are being investigated that can exploit low precision or mixed precision arithmetic.  In particular, these approaches use lower precision for storage and computation and benefit from distributed memory implementations.  New methods are being developed that can reduce global communication points, such as the computation of inner products. In another line of research, randomized numerical linear algebra has emerged as a powerful framework to drastically reduce computational costs.  Although initial randomized approaches that were based on obtaining low-rank approximations are not suitable for practical large-scale inverse problems, iterative methods that can exploit randomization for solving least-squares problems are being widely adopted and investigated.  

Inner-product-free Krylov subspace methods have recently emerged as competitive alternatives to well-known algorithms for inverse problems. In particular, the changing minimum residual Hessenberg (CMRH) method is a Krylov subspace method similar to generalized minimum residual (GMRES), but instead of using the Arnoldi scheme (in the case of GMRES) to construct an orthonormal basis of the Krylov subspace, CMRH uses a Hessenberg method to construct a good (although not orthonormal) linearly independent basis. Although the Hessenberg method was mentioned in the literature as early as 1950 (see, e.g., \cite[section 44]{faddeev1963}) for computing eigenvalues, it was Sadok \cite{sadok1999new} who introduced the CMRH method in 1999 as a way to solve linear systems.
In 2012, Sadok and Szyld \cite{sadok2012new} investigated the relationship between CMRH and GMRES. More recently, it was shown that CMRH is a regularizing iterative method, and that further regularization can be incorporated with a hybrid approach, called HCMRH \cite{brown2024hcmrh}. An extension to rectangular systems, called LSLU 
was developed in \cite{brown2024hlslu} that uses a generalized Hessenberg iterative algorithm to generate linearly independent bases for Krylov subspaces associated with $A^TA$ and $AA^T$. The approach has a theoretical connection to the LU factorization (with partial pivoting), and the iterative regularization properties of LSLU are very similar to LSQR \cite{paige1982lsqr,paige1982algorithm}. The main benefit of LSLU is that no inner-products are required during the iterations (note that we assume that matrix-vector mulitplications with $A$ and $A^T$ can be performed via function evaluations).  This leads to projected quasi-minimum residual problems that are smaller and easier to solve at each iteration. 

In general, a quasi-minimum residual method can provide suitable solutions and various relationships can be made between iterates \cite{saad2003iterative}.  However, in the context of inverse problems, the residual norm provides important information regarding the fit-to-data term, which depends on statistical assumptions about the measurement noise.  In particular, the solution to LS problem \eqref{eq:LS} corresponds to a maximum likelihood estimate and the solution to Tikhonov problem \eqref{eq:LS2} corresponds to a maximum a posteriori estimate \cite{calvetti2007introduction}.  Thus, we seek inner-product free Krylov methods that are also minimum residual methods. Unfortunately, this means that a tall, skinny least-squares problem needs to be solved at each iteration.  We address this computational challenge by exploiting recent work on randomized methods, which produces solutions that can approximate closely the minimal residual norm at each iteration.

\paragraph{Main Contributions}
We propose a new family of 
Krylov subspace methods that are inherently inner-product free and produce solutions with a smaller residual norm than existing inner-product free Krylov subspace methods.  
These are based on a (generalized) Hessenberg method with partial pivoting, where either one or two sets of linearly independent vectors are constructed to span Krylov subspaces. Contrary to standard iterative methods, the new approach avoids reorthogonalization and does not require inner-products. Contrary to other inner-product free Krylov methods, we use randomized sketching to obtain a more accurate approximation of the objective function to be minimized: the norm of the projection of either the residual or the Tikhonov objective function in a Krylov subspace of increasing dimension. The potential implications of this for inverse problems is that for existing inner-product free Krylov methods, it may be challenging to select regularization parameters in a hybrid approach and the bounds for the residuals depend on the conditioning of the basis vectors.

\paragraph{Paper structure}
The paper is organized as follows. In Section \ref{sec:SS_LS} we give an overview of sketch-and-solve methods for least-squares problems, including choices of sketching matrices. In Section \ref{sec:Hessenberg}, we recall the Hessenberg method to build a linearly independent basis for different choices of Krylov subspaces without requiring inner-products. We present new inner-product free Krylov methods, sketched CMRH (sCMRH) and sketched LSLU (sLSLU), to solve the least-squares problem and introduce an adaptation that considers Tikhonov regularization. Finally, in Section \ref{sec:numerics}, we illustrate the effectiveness of sCMRH and sLSLU, compared to existing iterative methods, and of sLSLU with Tikhonov regularization using a fixed regularization parameter. We provide closing remarks in Section \ref{sec:conc}\section{Background on randomized methods for least-squares problems}\label{sec:SS_LS}
In Section \ref{subsec:sLSLU} we use randomized methods to efficiently solve strongly overdetermined projected LS problems, so the aim of this section is to give some background on randomized methods.
Randomized numerical linear algebra, particularly those approaches involving sketching, has gained increasing popularity, see e.g \cite{Martinsson_Tropp_2020}. Sketching is a linear dimensionality reduction technique, and there are different ways in which sketching has been used to solve least squares problems. One of the conceptually simplest approaches is to sketch-and-solve, which can be used to find approximate solutions of least-squares problems where the system matrix is tall and skinny. This was originally proposed in \cite{4031351} and has gained a lot of attention due to its simplicity and probabilistic guarantees. In particular, we can define a sketching matrix $S \in \mathbb{R}^{\ell \times m}$, such that the following subspace (oblivious) embedding property is satisfied for any vector $a\in\mathbb{R}^{m}$ in a given set of vectors
\begin{equation}\label{eq:sub_eq}
    (1-\epsilon) \|a\| \leq \|Sa\|\leq (1+\epsilon)\|a\|.
\end{equation}
For this to be a dimensionality reduction technique, we typically assume that $\ell \ll m$. Even if this is a very favorable property, it is not trivial to construct such matrices deterministically in practice, in the sense that \eqref{eq:sub_eq} is guaranteed for any $a\in\mathbb{R}^{m}$. However, the analytical properties of random matrices can be used to construct sketching matrices $S \in \mathbb{R}^{\ell \times m}$ that will satisfy \eqref{eq:sub_eq} with high probability. Note that this is a special case of a random subspace embedding; for a formal definition, see, e.g. \cite[Chapter 8.1]{Martinsson_Tropp_2020}. 

Moreover, there exist different choices of random matrices in the literature that can be computationally cheap to construct and apply. In particular, we only consider the embedding of oblivious subspaces, which do not assume any prior information about the set of possible vectors $a$. The easiest class to analyze and implement is that of Gaussian embeddings, where each entry of $S \in \mathbb{R}^{\ell \times m}$ is an independent draw of a Gaussian distribution with zero mean and variance $1/\ell$. Note that applying a Gaussian sketch to a vector has an  $O(\ell m)$ cost (the explicit storage of the sketch matrix is also $O(\ell m)$), see \cite[Chapter 8.3]{Martinsson_Tropp_2020}. When dealing with high dimensional problems, one can also use structured random embeddings, which can reduce the storage and application costs. The most well-used sketches in this case are the subsampled randomized trigonometric transforms (SRTT), subsampled random Fourier transform (SRFT), and the sparse sign embeddings. However, the latter sketches require sufficiently efficient implementations to be faster than a simple Gaussian sketch in practice. 
For simplicity, in this paper we will only use Gaussian embeddings, but all results can be generalized to the use of other sketching techniques. 

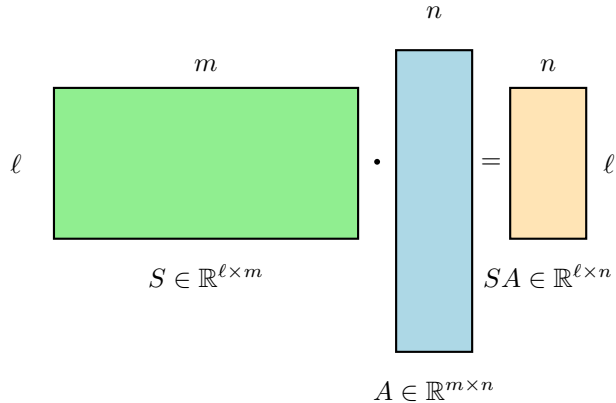
\begin{figure}[ht]
\centering
\begin{tikzpicture}[>=latex]
\definecolor{matrixA}{RGB}{173,216,230} 
\definecolor{matrixPhi}{RGB}{144,238,144} 
\definecolor{matrixB}{RGB}{255,228,181} 
\fill[matrixPhi] (0,1.5) rectangle (4,3.5); 
\draw[thick] (0,1.5) rectangle (4,3.5); 
\node at (2,1) {$S \in \mathbb{R}^{\ell \times m}$};
\fill[matrixA] (4.5,0) rectangle (5.5,4); 
\draw[thick] (4.5,0) rectangle (5.5,4); 
\node at (5,-0.5) {${A} \in \mathbb{R}^{m \times n}$};
\fill[matrixB] (6,1.5) rectangle (7,3.5); 
\draw[thick] (6,1.5) rectangle (7,3.5); 
\node at (6.5,1) {${SA} \in \mathbb{R}^{\ell \times n}$}; 
\node at (4.25,2.5) {\huge$\cdot$}; 
\node at (5.75,2.5) {$=$};     
\node at (-0.5,2.5) {$\ell$}; 
\node at (2,3.8) {$m$};    
\node at (5,4.5) {$n$};  
\node at (6.5,3.8) {$n$};  
\node at (7.3,2.5) {$\ell$};  
\end{tikzpicture}
\caption{ Schematic representation of the sketching of a matrix $A$ using a sketch $S$.}\label{fig:sketching}
\end{figure}

One of the most popular uses of sketching is to find approximate solutions of least-squares problems. Suppose we have a tall and skinny matrix $A \in \mathbb{R}^{m \times n}$, where $m \gg n$, then we can define a sketching matrix $S \in \mathbb{R}^{\ell \times m}$, such that $\ell \ll m$ and is assumed to be a small multiple of $n$. A representation of the sketched matrix $SA \in \mathbb{R}^{\ell \times n}$, of smaller dimension than $A$, can be observed in Figure \ref{fig:sketching}.
The sketch-and-solve method to find an approximate solution to the least-squares problem \eqref{eq:LS} involves solving the following minimization problem
\begin{equation} \label{sketched_LS}
\min_{x \in \mathbb{R}^n}  \| S(A x - b)\|, 
\end{equation}
where $S$ is a sketch matrix, see, e.g. \cite{4031351}, \cite[Chapter 10.3]{Martinsson_Tropp_2020}. Note that, when using Gaussian sketching, the solution of \eqref{sketched_LS} is an unbiased estimator of the solution of the original least-squares problem \eqref{eq:LS} provided that the matrix $A$ has full column rank. 

The fact that $A$ needs to be tall and skinny seems to be a restrictive property, since a lot of applications do not give rise to systems where the matrices are naturally in this form.  However, it has been proposed to use randomized techniques in combination with Krylov methods, where the projections give rise to such tall and skinny matrices. 
Specifically, for inverse problems, we know that a good approximation of the solution can be found in a Krylov subspace of small dimension involving the right-hand side $b$, $A$, and possibly $A^T$. 
Thus, as we will describe in Section \ref{sec:Hessenberg}, one can (1) construct a basis for the relevant Krylov subspace(s) and (2) use sketching to solve the projected least-squares problem, since this now involves a tall and skinny matrix even if $A$ was not. The Krylov basis need not be orthogonal, although it should be moderately well-conditioned, so the question is how to construct such a nonorthogonal Krylov basis. Previous approaches have considered using partial reorthogonalization, i.e. truncated (or incomplete) Golub-Kahan or Arnoldi methods, e.g. \cite{Guttel2023sketch}. However, these require setting an extra hyperparameter, namely the number $j$ of vectors one wants to reorthogonalize against (note that $j$ is usually taken to be $O(1)$).  Moreover, inner-products might still be related to algorithmic break-downs in low precision arithmetic due to information loss in the orthogonal projections, see \cite{brown2024hlslu} for details. For these reasons, in this paper we propose to use a completely inner-product free construction. 

It is worth mentioning that sketching has been used in other contexts as well.  Although it is beyond the scope of this paper to give a complete survey of methods, we mention a few examples that relate sketching with iterative Krylov methods. First, the approach of subspace embeddings has been used to sketch the inner products in the Gram-Schmidt orthogonalization procedure \cite{BalabanovRGMRES2022}. Second, a lot of research has also been done in the use of sketching as preconditioning for iterative methods, for example  in the LSNR \cite{lsnr} and Blendenpik \cite{Avron2010Blendenpik} algorithms. Both approaches can be more expensive than the sketch-to-solve approach, and target scenarios where high accuracy is required. Since this is typically not the case for inverse problems where the presence of noise and ill-conditioning limit reconstruction accuracy, neither of these paradigms will be addressed in this paper.

In fact, recently, some works have used sketching to do (almost orthogonal) projections onto Krylov subspaces, see, e.g. \cite{NakatsukasaSGMRES2024}. In the following section, we give an alternative to construct the basis for the Krylov subspaces that is inherently inner-product free and based on the Hessenberg method with (partial) pivoting.  Note that it is well accepted that partial pivoting is stable in practice, even so it is not backwards stable and unstable in the worst case; this is deemed to be very unlikely. In other words, we refer to the following quote attributed to Wilkinson: \textit{“Anyone that unlucky has already been run over by a bus!”}.

\section{Inner-product free Krylov subspace methods}\label{sec:Hessenberg}
Krylov iterative algorithms are a class of very powerful projection methods that make use of Krylov subspaces of the form:
$$ {\cal{K}}_{k}(C,d) = span \{d, Cd,C^2d, ...,C^{k-1}d \},$$ 
for a given matrix $C$ and a vector $d$. In their most general form, they can be classified by the (Krylov) solution subspace, where a solution at each iteration is sought. The optimality conditions that are imposed determine the solution in that subspace. For orthonormal bases, the optimality conditions are usually given in terms of the orthogonality of the residual with respect to another Krylov subspace.

In this section, we focus on inner-product free Krylov methods, where the first step is to build a nonorthogonal basis for the relevant Krylov subspaces. In \Cref{sub:LSLU} we describe CMRH and LSLU, two inner-product free quasi-minimum residual methods.  A key component of these approaches is the (generalized) Hessenberg method used to generate bases for Krylov subspaces without using inner-products and in a practically stable manner. Then, in \Cref{subsec:sLSLU}, we describe our proposed approaches sCMRH and sLSLU, which are inner-product free methods that seek to minimize the residual norm or Tikhonov objective function, where sketching is used to solve the projected problem efficiently.

\subsection{CMRH and LSLU: quasi-minimum residual methods}\label{sub:LSLU}
As described in the introduction, the Hessenberg method is an inner-product free approach to generate a basis for ${\cal{K}}_{k}(A,b)$ \cite{sadok1999new,sadok2012new} where $A$ is a square matrix.  That is, 
after $k$ iterations, we have a unit lower triangular matrix $L^{(s)}_{k+1}$ and an upper Hessenberg matrix $H^{(s)}_{k+1,k}$ that satisfy
\begin{eqnarray}\label{Hessenberg_s}
    AL^{(s)}_k &=& L^{(s)}_{k+1}H^{(s)}_{k+1,k},
\end{eqnarray}
where the range of $L^{(s)}_k$ is ${\cal{K}}_{k}(A,b)$. 
For rectangular matrices $A$, the generalized Hessenberg process \cite{brown2024hlslu} was established as an inner-product free method to generate bases for both ${\cal{K}}_{k}(A^T A,A^T b)$ and ${\cal{K}}_{k}(A A^T,b)$. In this case, the following two relations are updated at each iteration,
\begin{eqnarray}\label{Hessenberg_r}
    AL^{(r)}_k &=& D_{k+1}H^{(r)}_{k+1,k} \nonumber \\ 
    A^TD_{k+1} &=& L^{(r)}_{k+1}W_{k+1},
\end{eqnarray}
where the range of $L^{(r)}_k$ is ${\cal{K}}_{k}(A^T A,A^T b)$, the range of $D_{k+1}$ is ${\cal{K}}_{k}(A A^T,b)$ and both matrices are unit lower triangular. Moreover, $H^{(r)}_{k+1,k}$ is upper Hessenberg and $W_{k+1}$ is upper triangular. For the analysis of this method, see \cite{brown2024hlslu}.

In the previous paragraph, we used the superscript $^{(s)}$ for quantities resulting from applying the Hessenberg method to {\em square} matrices, and the superscript $^{(r)}$ for quantities resulting from applying the generalized Hessenberg method to {\em rectangular} matrices. But from now on we will drop the superscript for $L_{k}$ (and $H_{k+1,k}$), under the general assumption that the range of $L_k$ is a relevant Krylov subspace, and that both matrices are obtained after $k$ iterations of either the Hessenberg or the generalized Hessenberg method.

Finding the minimal residual norm solution in the corresponding Krylov subspace spanned by the columns of $L_{k}$ (and $x_0$) corresponds to solving:
\begin{equation} \label{original_LS_proj}
y^{(LS)}_k = \arg\min_{y \in \mathbb{R}^k}  \| A L_{k}y - b\|, \quad x^{(LS)}_k = L_{k}y^{(LS)}_k +x_0.
\end{equation}
However, the columns of $L_k$ (and, if relevant, $D_k$) are not orthonormal. Therefore, it is not straightforward to use the relations in \eqref{Hessenberg_s} or \eqref{Hessenberg_r} to reduce the least-squares problem \eqref{original_LS_proj} into a small projected problem. 
This is different from other standard Krylov methods, like GMRES and LSQR, where the Arnoldi relationship or the Golub-Kahan bidiagolization can be used to solve \eqref{original_LS_proj} efficiently \cite{golub2013matrix,saad2003iterative}. Note that this is the main drawback of this type of Krylov basis construction. However, recall that the basis vectors are computed in a stable way in practice using (partial) pivoting.  Moreover, the approach is inherently inner-product free.

There are Krylov subspace methods that use the relations \eqref{Hessenberg_s} or \eqref{Hessenberg_r} to approximate the orthogonal projection \eqref{original_LS_proj} by an oblique projection of the residual. For example, the CMRH method for square matrices $A$ \cite{brown2024hcmrh,sadok1999new,sadok2012new}, and the LSLU method for rectangular matrices \cite{brown2024hlslu}
seek to
approximate the least squares solution by computing,
\begin{equation} \label{Hessenber_LS}
x_k =\arg\min_{x \in x_0+R(L_k)}  \| M_{k+1}^{\dagger} (Ax-b) \|,
\end{equation}
where $M_{k+1}=L^{(s)}_{k+1}$ for CMRH and $M_{k+1}=D_{k+1}$ for LSLU, and $M^{\dagger}_{k+1}$ is the Moore-Penrose pseudoinverse of $M_{k+1}$. The residual norm associated with this solution $x_k$ fulfills the following inequality:
\begin{equation}\label{eq:eqthm} 
\| Ax^{(LS)}-b \| \leq \|A x_k - b \| \leq \kappa(M_{k+1})\, \| Ax^{(LS)}-b \|  \end{equation}
where $\kappa(\cdot)$ denotes the condition number and from \eqref{original_LS_proj},
$x^{(LS)}$ is the minimal norm solution in the given Krylov subspace: $\|A x^{(LS)} - b \| = \min_{x \in x_0+R(L_k)}  \| Ax-b \| $. For more details, see \cite{brown2024hcmrh} for CMRH and \cite{brown2024hlslu} for LSLU.

Note that an interesting interpretation of CMRH and LSLU, which was not noted in the original papers, is that it corresponds to the following approximation of \eqref{original_LS_proj},
$$x^{(LS)}_k=L_k y^{(LS)}_k = L_k (M_{k+1} H_{k+1,k})^{\dagger}b \approx L_k H_{k+1,k}^{\dagger} M_{k+1}^{\dagger}b.$$
In other words, an approximate solution for \eqref{original_LS_proj} can be obtained in two stages,
\begin{eqnarray}\label{LSLU_alt}
    z_k & =& \arg\min_{z \in \mathbb{R}^{k}}  \| M_{k+1} z - b\|, \\
    x_k & =& L_k (H_{k+1,k}^T H_{k+1,k})^{-1} H_{k+1,k}^T z_k.
\end{eqnarray} 
However, the approximation of the pseudoinverse of the product as  $(M_{k+1} H_{k+1,k})^\dagger \approx H_{k+1,k}^{\dagger} M_{k+1}^{\dagger},$ is not true in general.
Assuming no breakdown, $M_{k+1}$ has full column rank, so the approximation is exact in the case where $H_{k+1,k}$ has full row rank.  In the next section, we consider a sketched approximation of $(M_{k+1} H_{k+1,k})^{\dagger}.$

\subsection{sCMRH and sLSLU: sketched inner-product free methods}
\label{subsec:sLSLU}
In this paper, we propose new methods that use the (non-orthogonal) basis vectors for the relevant Krylov subspace, generated with the Hessenberg or the generalized Hessenberg method,
and approximately compute the solution of the projected least squares problem using a sketch-and-solve approach.  The sketched, projected problem is given by
\begin{equation} \label{sketched_LSLU}
\min_{y \in \mathbb{R}^k}  \| S(A L_{k}y - b)\|, 
\end{equation}
where $S$ is an appropriate sketching matrix.  Assuming $x_0 = 0$, the solution at the $k$th iteration is given by $x^{(S)}_k=L_k y^{(S)}_k$, where $y^{(S)}_k$ is the solution of \eqref{sketched_LSLU}. The crux of our approach lies in the assumption that at any given iteration, the solution $x^{(S)}_k$ will be close to the minimal residual norm solution $x^{(LS)}_k$ with high probability. 

Using the subspace embedding property \eqref{eq:sub_eq} on the residual norm, and recalling that 
$x^{(LS)}$ is the minimal norm solution in the given Krylov subspace, i.e., $$\|A x^{(LS)} - b \| = \min_{x \in x_0+R(L_k)}  \| Ax-b \|,$$
then one can easily see that
\begin{equation}\label{eq:eqthm2} 
\| Ax^{(LS)}-b \| \leq \|A x^{(S)}_k - b \| \leq \frac{1+\epsilon}{1-\epsilon} \, \| Ax^{(LS)}-b \|.
\end{equation}
For example, if $S \in \mathbb{R}^{\ell\times m}$ is a Gaussian sketch, \eqref{eq:eqthm2} holds with a small probability of failure if $\ell \sim m \log(m)/\epsilon^2$ (in practice, this is usually further reduced to $\ell \sim m /\epsilon^2$).  Moreover, as noted in \cite{Meier2024Stable},  $\epsilon$ is the subspace embedding
constant of the matrix $[A\,\, b\,]$ (i.e. appending $b$ to the matrix $A$). This satisfies $\frac{1+\epsilon}{1-\epsilon} = \kappa_2 (S Q^{[A \, b]})$, where $Q^{[A \, b]}$ is the orthogonal matrix obtained doing the QR factorization of $[A\,\, b\,]$, and $\epsilon$ is usually of the order of $0.5$. An accurate mathematical description of the relevant theory can be found in \cite[Section 8.7.]{Martinsson_Tropp_2020}, and sharper bounds with a different probability of failure can be found in the seminal paper \cite{2006SarlosSketch}. 

This means that, in theory one can pick $\epsilon$ to be small enough that the solution $x_k^{(S)}$ produced by either sCMRH and sLSLU will have a smaller residual norm than that of the solution $x_k$ obtained using the inner-product free (but not sketched) counterparts. Although this is a very strong motivation of this work, it cannot, however, always be guaranteed in practice.

Moreover, we know that if $S$ is a Gaussian sketch, then the sketch and solve solution to the projected problem is an unbiased estimate for the least-squares solution.
Consider the projected problem and the sketched, projected problem,
$$\min_{y \in \mathbb{R}^k}  \| A L_{k}y - b\| \quad \mbox{and}\quad \min_{y \in \mathbb{R}^k}  \| S(A L_{k}y - b)\|,$$
respectively, where $A L_k$ is a full column rank matrix.
Then
$$\mathbb{E}[(S A L_{k})^\dagger (S b)] = (A L_{k})^\dagger b$$
where $\dagger$ denotes the Moore-Penrose pseudoinverse. Since $L_k$ is independent of the sketch, we have that $\mathbb{E}[ x_k^{(S)}] = \mathbb{E}[ L_k (S A L_{k})^\dagger (S b)] = L_k \mathbb{E}[(S A L_{k})^\dagger (S b)]= L_k (A L_{k})^\dagger b = x_k^{(LS)}.$
In addition, we can further analyze the expected squared residual norm as 
\begin{align*}
\mathbb{E}[\| A L_k y_k^{(S)} - b\|^2] & = \left(1 + \frac{k}{\ell - k -1}\right) \min_y \|A L_k y - b \|^2 \\
& = (1+\varepsilon) \min_y \|A L_k y - b \|^2
\end{align*}
when the sketch dimension is $\ell = \frac{k}{\varepsilon} + k + 1.$ See \cite{derezinski2024recent,drineas2011faster} and references therein.

Algorithm \ref{alg:sLSLU} provides the details for sketched LSLU (sLSLU); the sketched version of CMRH can be straightforwardly generalized. Note that for sLSLU, one sketching matrix is needed $S \in \mathbb{R}^{\ell \times m}$; however, to incorporate Tikhonov regularization, we will need two sketching matrices, $S_1 \in \mathbb{R}^{\ell \times n}$ and $S_2 \in \mathbb{R}^{\ell \times m}$. 

\begin{algorithm}[h]
\caption{sLSLU} \label{alg:sLSLU}
\begin{algorithmic}[1]
\REQUIRE $A$, $b$, $x_0$, $\text{maxiter}$, $S_1$, $S_2$
\STATE Define $t = [1,2,\ldots,m]^T$, $g = [1,\ldots,n]^T $.
\STATE $r_0 = b - A x_0 $ 
\STATE Determine $i$ such that $|r_0(i)| = \|r_0\|_{\infty}$
\STATE $\beta = r_0(i)$; $d_1 = r_0/ \beta$; $t(1) \Leftrightarrow t(i)$
\STATE $v_0 = A^T r_0$
\STATE Determine $i_2$ such that $|v_0(i_2)| = \|v_0\|_{\infty}$
\STATE $\alpha = v_0(i_2)$; $l_1 = v_0/ \alpha$; $g(1) \Leftrightarrow g(i_2)$
\STATE $r = A^Td_1$; $W(1,1) = r(g(1))$
\FOR{$k = 1,\ldots,\text{maxiter}$} 
  \STATE $u = A^Td_{k}$
   \STATE $z_k = S_2 * u$
  \FOR{$j = 1,\ldots, k$}
  \STATE $H(j,k) = u(t(j))$ ; $u = u-H(j,k)d_j$
  \ENDFOR
  \IF{$k<m$ and $u \neq 0$}
  \STATE Determine $i \in \{ k+1,\ldots,m\}$ such that $|u(t(i))| = \|u(t(k+1:m))\|_{\infty}$
  \STATE $H(k+1,k) = u(t(i))$; $d_{k} = u/H(k+1,k)$; $t(k+1) \Leftrightarrow t(i)$

  \ELSE
   \STATE $H(k+1,k) = 0$;
  \ENDIF 
  \STATE $q = A^Td_k$
\FOR{$j = 1, \ldots,k$}
  \STATE $W(j,k+1) = q(g(j))$; $q = q-W(j,k+1)l_j$
  \ENDFOR
  \IF{$k<n$ and $q \neq 0$}
  \STATE Determine $i_2 \in \{ k+1,\ldots,n\}$ such that $|q(g(i_2))| = \|q(g(k+1:n))\|_{\infty}$
  \STATE $W(k+1,k+1) = u(g(i_2))$; $l_{k+1} = q/W(k+1,k+1)$; $g(k+1) \Leftrightarrow g(i_2)$
  \ELSE
 \STATE break
  \ENDIF   
 \STATE Compute $y_{k}$ to be the minimizer of $\| S_2r_0 - S_2AL_ky \|_2^2 = \| S_2r_0 - Zy \|_2^2 $
  \STATE $ x_k = x_0 + L_k y_{k}$
\ENDFOR 
\end{algorithmic} 
\end{algorithm}

As can be observed in Algorithm \ref{alg:sLSLU}, the generalized Hessenberg method with partial pivoting requires locating the largest absolute value of $r_0$, $v_0$, and two other vectors at each iteration, which can be costly due to global communications. To avoid this occurrence, a pivoting alternative was proposed in \cite{brown2024hlslu} where only a small subset of elements in those vectors is observed, and the pivot is taken to be the element with the largest absolute value from this subset. In practice, this leads to a reduced communication cost and a stable way of building the non-orthogonal basis for the relevant Krylov Subspace methods: in our experiments, we never observed the condition number of the basis to grow significantly large. However, this strategy might still increase the condition number $\kappa(M_{k+1})$ in equation \eqref{eq:eqthm}, rendering the solutions from LSLU and CMRH less accurate. Therefore, this is a particular case of LSLU and CMRH where the sketched methods have a great potential compared to their non-sketched counterparts. All numerical results will implement this technique for different sample sizes.

 From Theorem $2.1$ of \cite{brown2024hlslu}, the authors derived a bound on the difference between the residual norms of solutions computed using LSLU and LSQR. It was shown that if the condition number of $\hat{R}_{k+1}$, the upper triangular matrix from the QR decomposition of $D_{k+1}$, does not grow too quickly, then the residual norms associated with the approximate solutions of LSLU and LSQR at each iteration are close to each other. In \Cref{sec:numerics}, we will compare the residual norm of sLSLU with the residual norms from Theorem $2.1$ of \cite{brown2024hlslu}.

In the next subsection we propose a sketch-and-solve approach to project the Tikhonov problem on a Krylov subspace and approximately compute the solution of the least squares problem \eqref{sketched_LSLU}.

\subsection{Extensions for Tikhonov Regularization}
Consider the standard Tikhonov regularization problem \eqref{eq:LS2}. The LSLU method computes a solution at the $k$th iteration to the following optimization problem
\begin{equation} \label{eq:TP} 
\min_{x \in x_0+\mathcal{R}(L_k)} \| D^{\dagger}_{k+1}(b-Ax)\|^2 + \lambda^2 \| L^{\dagger}_{k} x \|^2, 
\end{equation}
where similar to LSLU, the residual norm is replaced by a semi-norm, and the regularization term also includes a semi-norm. This is equivalent to solving 
\begin{equation} \label{eq:TSP} 
y_{\lambda ,k} = \arg \min_{y \in \mathbb{R}^k} \| \beta e_1 - H_{k+1,k}y \|^2 +  \lambda^2 \|y\|^2, 
\end{equation} where $\beta$ is selected using the 
 following sample strategy in \cite{brown2024hlslu}: select a small random sample of entries from $r_0$, $v_0$, $u$, $q$ and choose the largest value (in magnitude) in that sample. Similarly, one can use a sketch-and-solve approach to project the Tikhonov problem on a Krylov subspace and approximately compute the solution of the least squares problem using the following expression:
\begin{equation} \label{sketched_LSLU_hybrid}
\min_{y \in \mathbb{R}^k}  \| S_2(A L_k y - b)\|+ \lambda^2 \| S_1 (L_{k} y)\|^2, 
\end{equation}
where $S_1$ and $S_2$ are  appropriate sketching matrices. 
Similar to LSQR and LSLU, the authors of \cite{brown2024hlslu} derived a bound on the difference between the residual norms of solutions computed using LSLU and LSQR for the Tikhonov problem. It was shown that if the condition number of $\bar{D}_{k+1}$ from the block matrix  \begin{equation}
\label{eq:Dbar}
\overline{D}_{k+1} = \begin{bmatrix} D_{k+1} & 0 \\ 0 & L_{k} \end{bmatrix} 
\end{equation} does not grow too quickly, then the residual norms associated to the solution of LSLU for the Tikhonov problem is close to the residual norm of the solution obtained with LSQR for the Tikhonov problem (See Theorem $3.1$ of \cite{brown2024hlslu}). To display behavior of sLSLU with Tikhonov regularization, we will plot the associated residual norm against the residuals norms denoted in Theorem $3.1$ of \cite{brown2024hlslu} (See \Cref{sec:numerics}). An implementation of sLSLU with Tikhonov regularization is provided in \Cref{alg:sHLSLU}, which corresponds to \Cref{alg:sLSLU} if $\lambda$ = 0. Note that for all numerical results, we select a fixed regularization parameter. The correspondent implementation of sCMRH with Tikhonov regularization can be easily generalized.

\begin{algorithm}[H]
\caption{sLSLU with Tikhonov Regularization} \label{alg:sHLSLU}
\begin{algorithmic}[1]
\REQUIRE $A$, $b$, $x_0$, $\text{maxiter}$, $S_1$, $S_2$, $\lambda$
\STATE Define $t = [1,2,\ldots,m]^T$, $g = [1,\ldots,n]^T $.
\STATE $r_0 = b - A x_0 $ 
\STATE Determine $i$ such that $|r_0(i)| = \|r_0\|_{\infty}$
\STATE $\beta = r_0(i)$; $d_1 = r_0/ \beta$; $t(1) \Leftrightarrow t(i)$
\STATE $v_0 = A^T r_0$
\STATE Determine $i_2$ such that $|v_0(i_2)| = \|v_0\|_{\infty}$
\STATE $\alpha = v_0(i_2)$; $l_1 = v_0/ \alpha$; $g(1) \Leftrightarrow g(i_2)$
\STATE $r = A^Td_1$; $f_1 = S_1 *r$; $W(1,1) = r(g(1))$
\FOR{$k = 1,\ldots,\text{maxiter}$} 
  \STATE $u = A^Td_{k}$
   \STATE $z_k = S_2 * u$
  \FOR{$j = 1,\ldots, k$}
  \STATE $H(j,k) = u(t(j))$ ; $u = u-H(j,k)d_j$
  \ENDFOR
  \IF{$k<m$ and $u \neq 0$}
  \STATE Determine $i \in \{ k+1,\ldots,m\}$ such that $|u(t(i))| = \|u(t(k+1:m))\|_{\infty}$
  \STATE $H(k+1,k) = u(t(i))$; $d_{k} = u/H(k+1,k)$; $t(k+1) \Leftrightarrow t(i)$

  \ELSE
   \STATE $H(k+1,k) = 0$;
  \ENDIF 
  \STATE $q = A^Td_k$
  \STATE $f_{k+1} = S_1q$
\FOR{$j = 1, \ldots,k$}
  \STATE $W(j,k+1) = q(g(j))$; $q = q-W(j,k+1)l_j$
  \ENDFOR
  \IF{$k<n$ and $q \neq 0$}
  \STATE Determine $i_2 \in \{ k+1,\ldots,n\}$ such that $|q(g(i_2))| = \|q(g(k+1:n))\|_{\infty}$
  \STATE $W(k+1,k+1) = u(g(i_2))$; $l_{k+1} = q/W(k+1,k+1)$; $g(k+1) \Leftrightarrow g(i_2)$
  \ELSE
 \STATE break
  \ENDIF   
 \STATE Compute $y_{\lambda_k,k}$ to be the minimizer of $\| S_2r_0 - S_2AL_ky \|_2^2 +  \lambda^2\|S_1L_ky\| $
  \STATE $ x_k = x_0 + L_k y_{\lambda_k,k}$
\ENDFOR 
\end{algorithmic} 
\end{algorithm}

\section{Numerical Results}\label{sec:numerics}
In this section, we illustrate the effectiveness of sketched inner-product-free Krylov methods:  sketched CMRH (sCMRH) and sketched LSLU (sLSLU) in comparison to their non-sketched counterparts CMRH and LSLU \cite{brown2024hlslu}, and the classical GMRES and LSQR. 
We utilize three different test problems: a deblurring problem, a neutron tomography simulation from the IR Tools package \cite{2019GazzolaIRtools}, and an example with real data consisting of two open access datasets from the Finnish Inverse Problems Society \cite{bubba2017tomographic,hamalainen2015tomographic}. Moreover, for Tikhonov regularization with fixed $\lambda$, we provide numerical results for sLSLU, LSLU, and LSQR for some of these test problems. Note that each sketch matrix has the following dimensions: $S_1 \in \mathbb{R}^{\ell \times n}$, $S_2 \in \mathbb{R}^{\ell \times m}$  where $\ell = 10*(max\; iteration +1)$.

\subsection{Deblurring problem}
The first experiment consists of a deblurring experiment, where the aim is to reconstruct MATLAB's test image `cameraman' of size 256 $\times$ 256 pixels, which was corrupted by motion blur and additive Gaussian noise. The forward model was simulated using IR Tools \cite{2019GazzolaIRtools}, and noise was added so that the observation contained a $1\%$ noise level. 

\begin{figure}[ht]
\centering
\begin{tabular}{cc}
    {\scriptsize Relative Reconstruction Error Norms } &  {\scriptsize Residual Norms} \\
     \includegraphics[width=0.45\textwidth]{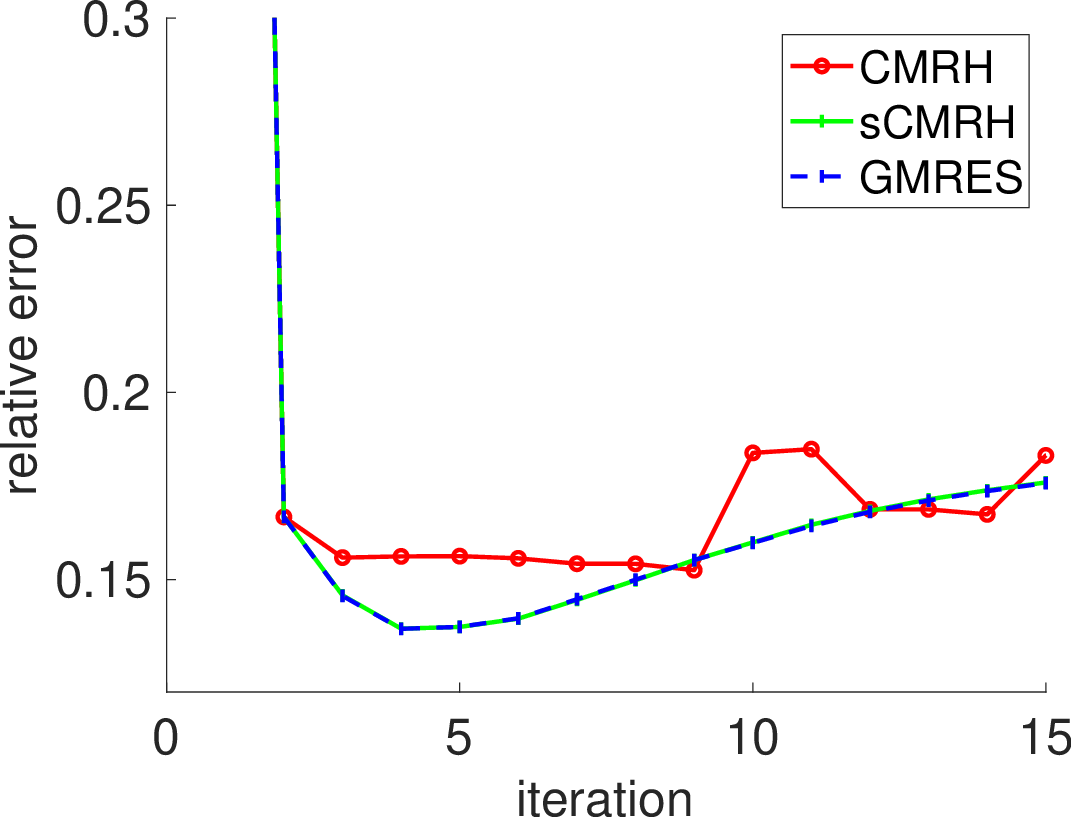} &  
   \includegraphics[width=0.45\textwidth] {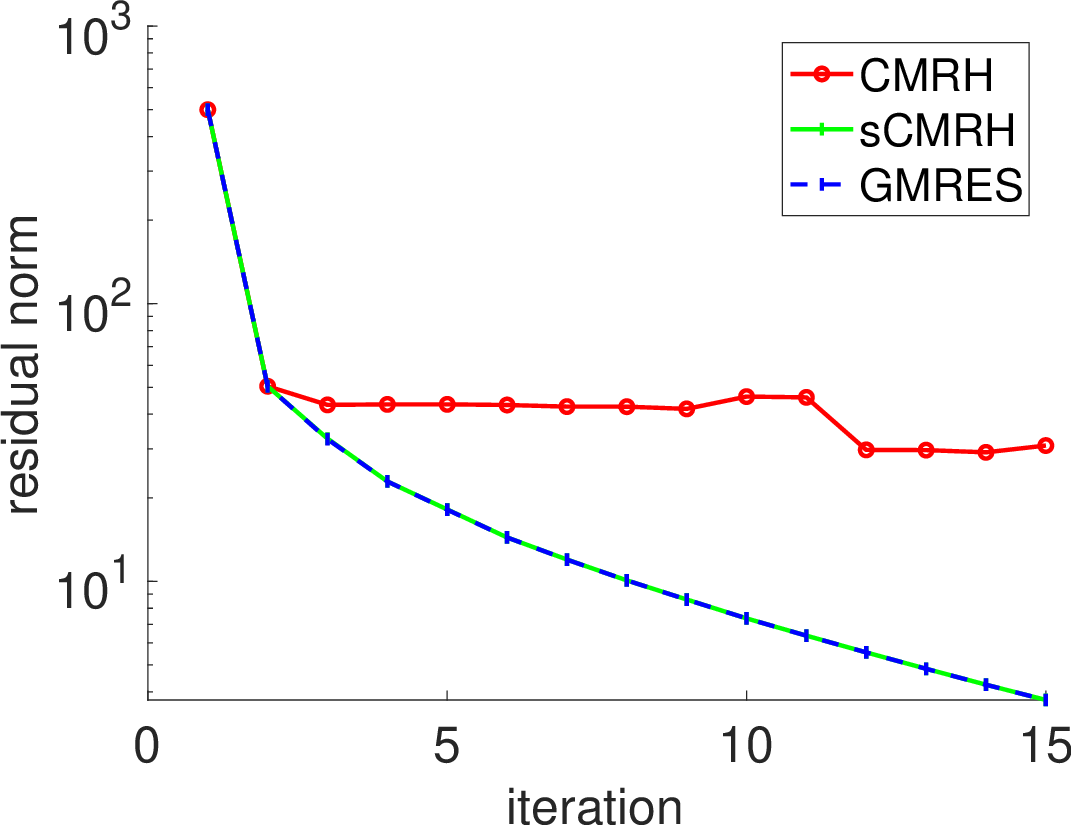} \\  
\end{tabular}
\caption{Relative reconstruction error norms (left) and residual norms (right).
In this case, sketched CMRH uses the pivots dictated by the maximum absolute value  from a set of randomly sampled coefficients ($5$).}
\label{fig:scmrh}
\end{figure}

\begin{figure}[ht]
\centering
\title{Cameraman}
\begin{tabular}{cccc}
    {\scriptsize Noisy Data} &
    {\scriptsize GMRES} &  {\scriptsize CMRH} & {\scriptsize sCMRH}\\ 
    \includegraphics[width=2.7cm]{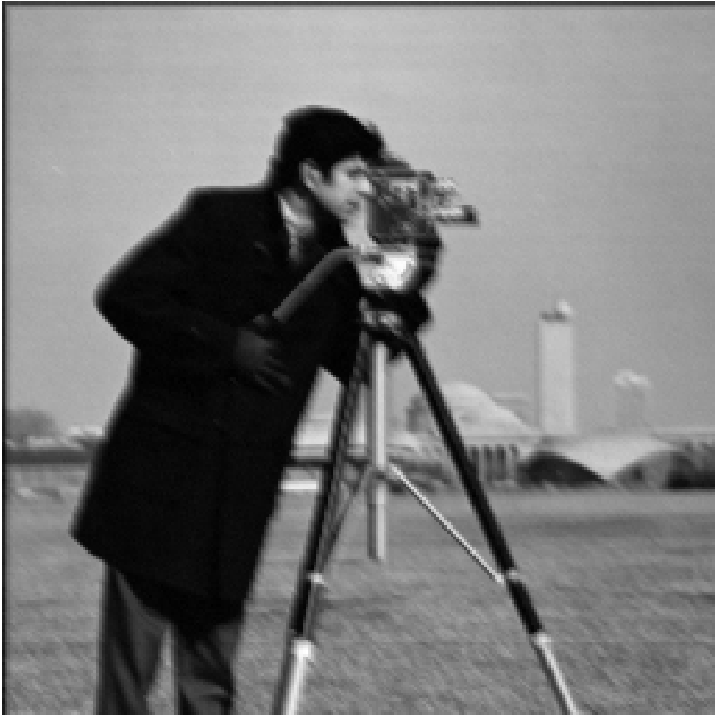} & 
    \includegraphics[width=2.7cm]{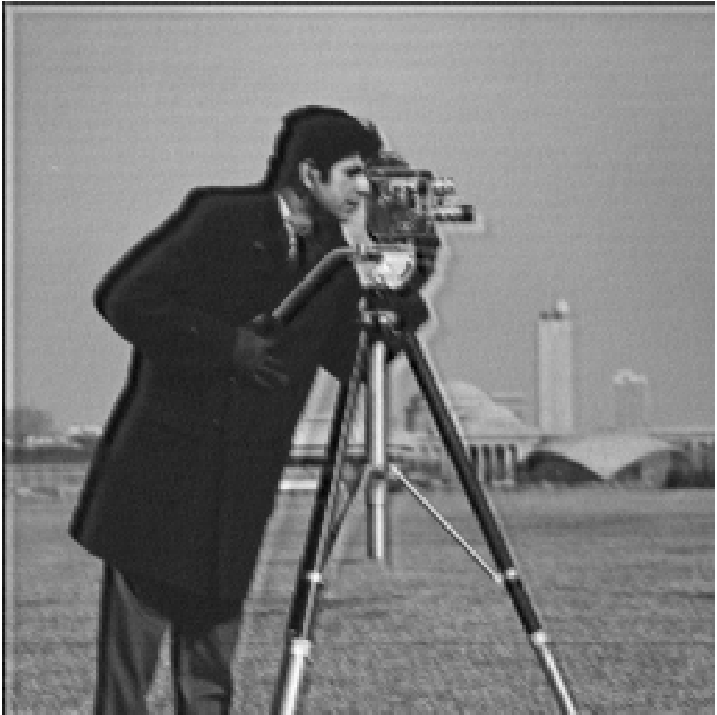} &  \includegraphics[width=2.7cm]{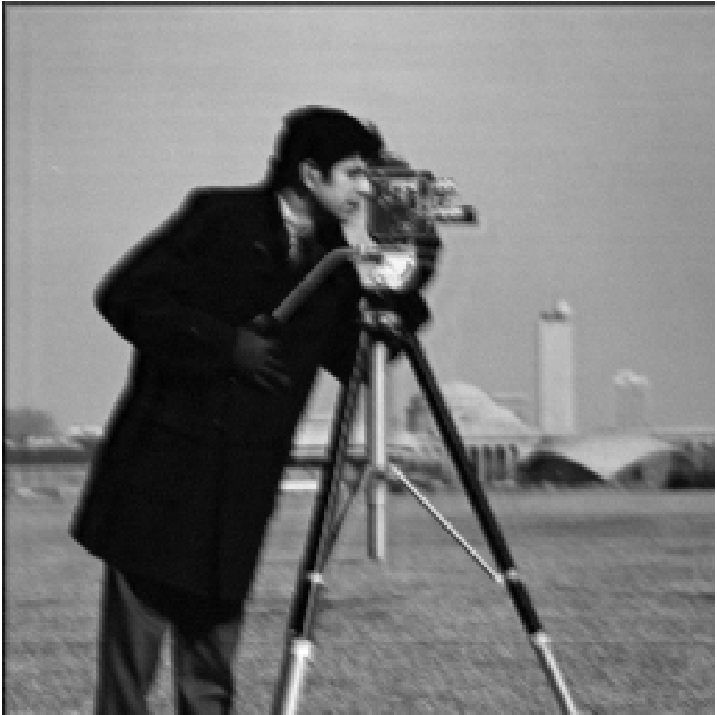} &
    \includegraphics[width=2.7cm]{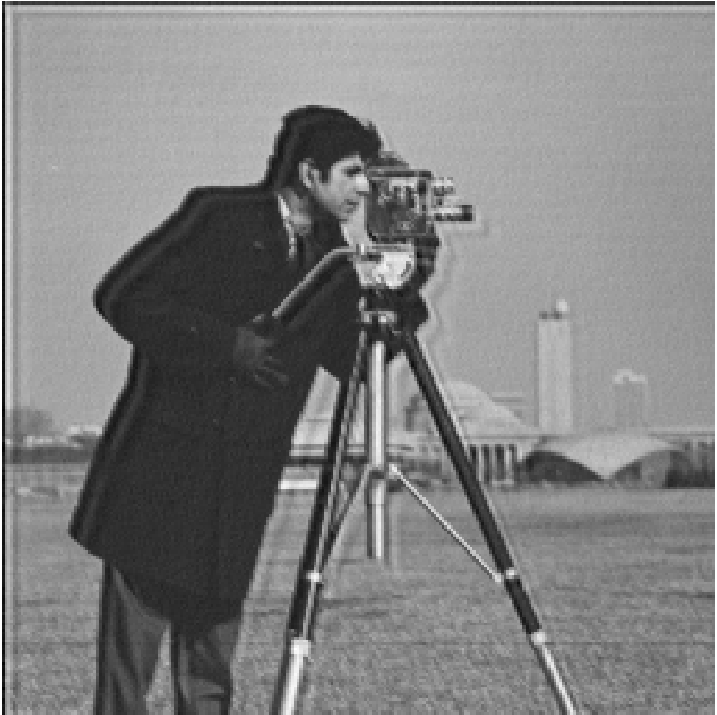} \\
\end{tabular} \\
\caption{Measured noisy data, and reconstructed images using GMRES, CMRH, and sCMRH.}
\label{fig:cameraman}
\end{figure}

Since this is a square problem, we compare sCMRH to CMRH and GMRES. The relative reconstruction error norm and residual norms per iteration are are displayed in Figure \ref{fig:scmrh}. We observe that the curves corresponding to sCMRH closely resembles that of GMRES, as dictated by the theory, while the residual norms for CMRH deviate as the iterations proceed. This can also be observed in the relative error norms. Thus, sCMRH produces solutions that more closely resemble GMRES solutions, but without the need for inner-product computations.  The reconstructed images are shown in Figure \ref{fig:cameraman}.

\subsection{Neutron tomography simulation}
Neutron imaging is a tomography technique based on gamma-rays that allows us to inspect the interior of dense or metallic objects. This is because, contrary to X-ray tomography, the absorption of neutrons is higher in `light' elements and lower in metallic elements. See, for example, the interior of a padlock in \cite{biguri2024tigrev3}. Note that, mathematically, this CT modality has the same mathematical model as X-ray tomography, but using different absorption coefficients for each material. Since most datasets for this type of tomography are proprietary, in this example, we use MATLAB's built-in demo image `circuit.tiff', which has similar structure to neutron tomography examples.

\begin{figure}[ht]
\centering
\begin{tabular}{cc}
    {\scriptsize Relative Reconstruction Error Norms } &  {\scriptsize Residual Norms} \\
     \includegraphics[width=0.45\textwidth]{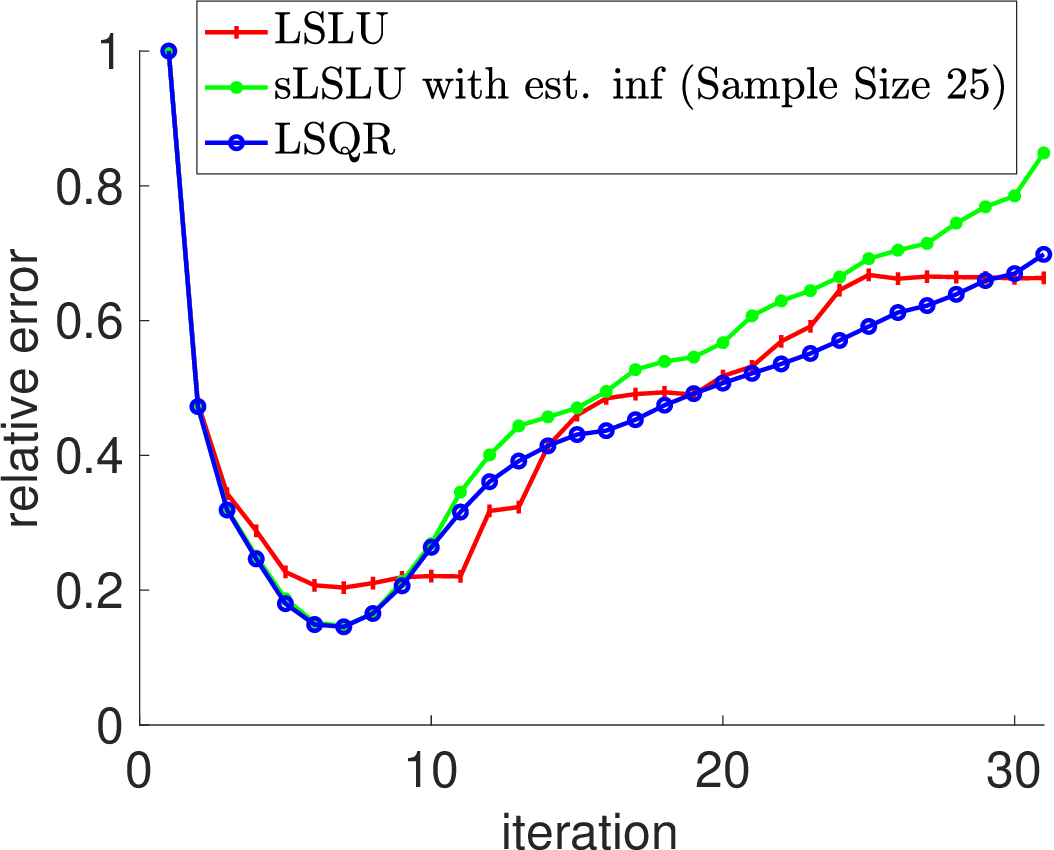} &  
   \includegraphics[width=0.45\textwidth] {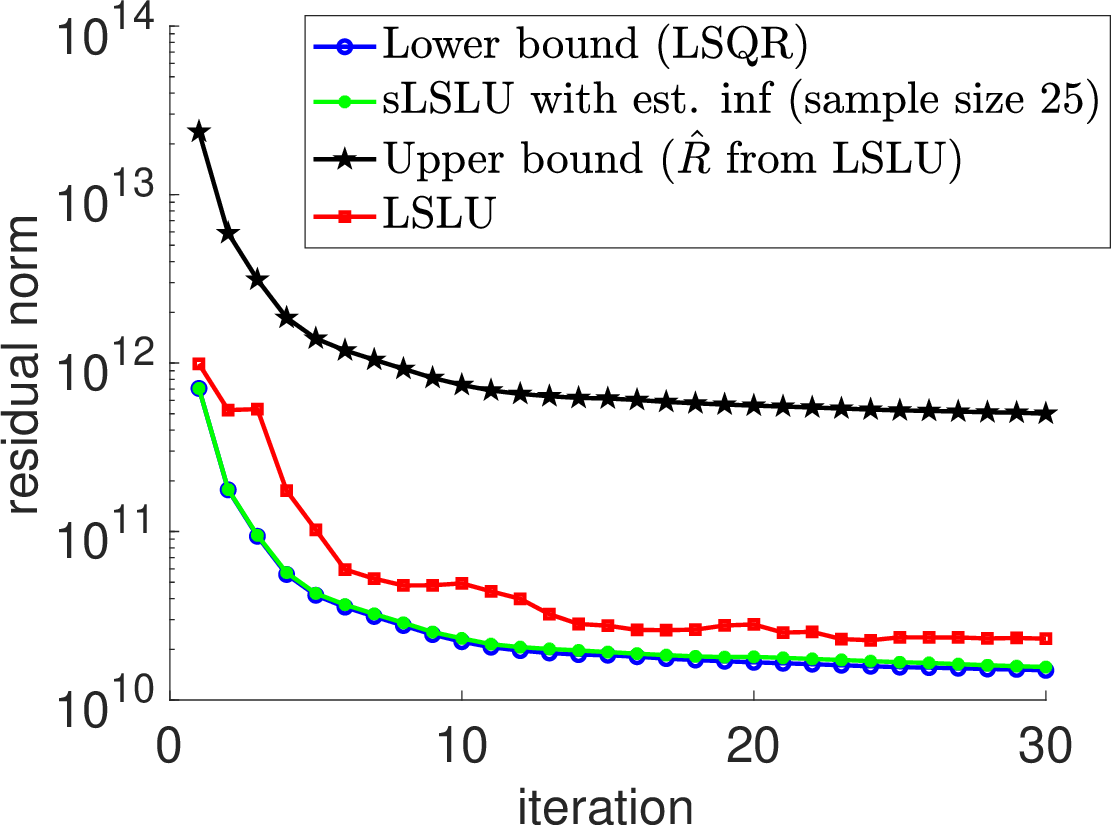}
\end{tabular}
\caption{Relative reconstruction error norms (left) and residual norms (right). In this case, sketched
LSLU uses the pivots dictated by the maximum absolute value from a set of randomly sampled
coefficients (25).}
\label{fig:sLSLU_rel}
\end{figure}

The sketched LSLU (sLSLU) algorithm implements the pivoting alternative from \cite{brown2024hlslu}. The sample size to approximate the infinity norm contains $25$ entries. Note that varying the sample size does not appear to drastically change the numerical results.
In \Cref{fig:sLSLU_rel}, we provide the relative reconstruction error norms per iteration of sLSLU. Results for LSQR and LSLU are provided for comparison. We observe that sLSLU performs better than LSLU and similar to LSQR, especially in the earlier iterations. Provided we implement a ``good" stopping criteria, we can compute an approximation that is of comparable quality to that produced with LSQR; while avoiding inner-products.

\begin{figure}[ht]
\centering
\title{Neutron Tomography Simulation} \\[12pt]
\begin{tabular}{cc}
    {\scriptsize Noisy Data} & {\scriptsize True Solution} \\
    \includegraphics[width=3.0cm]{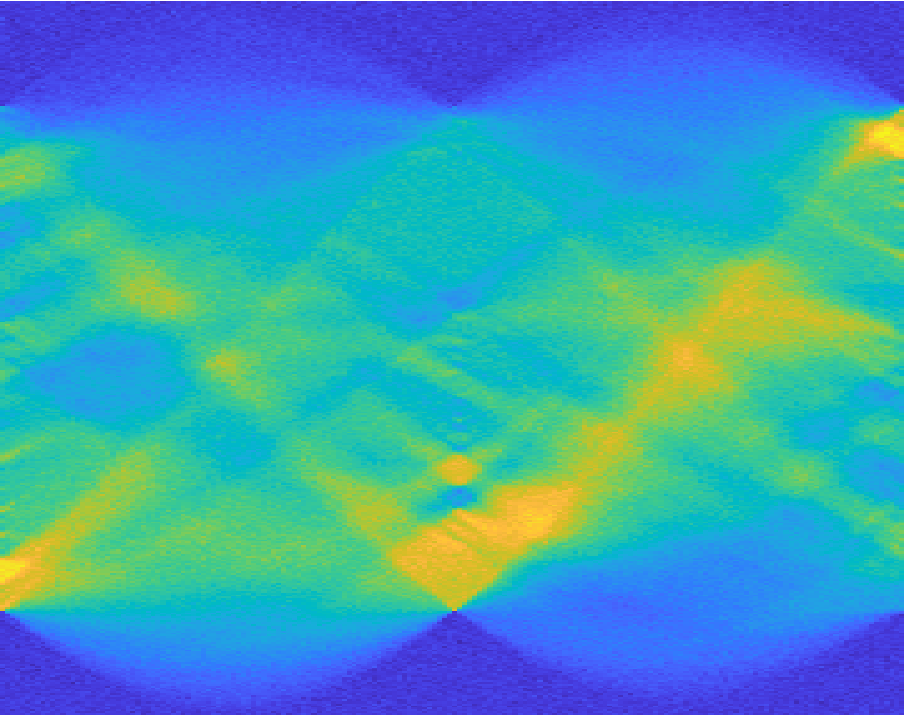} &
    \includegraphics[width=3.0cm]{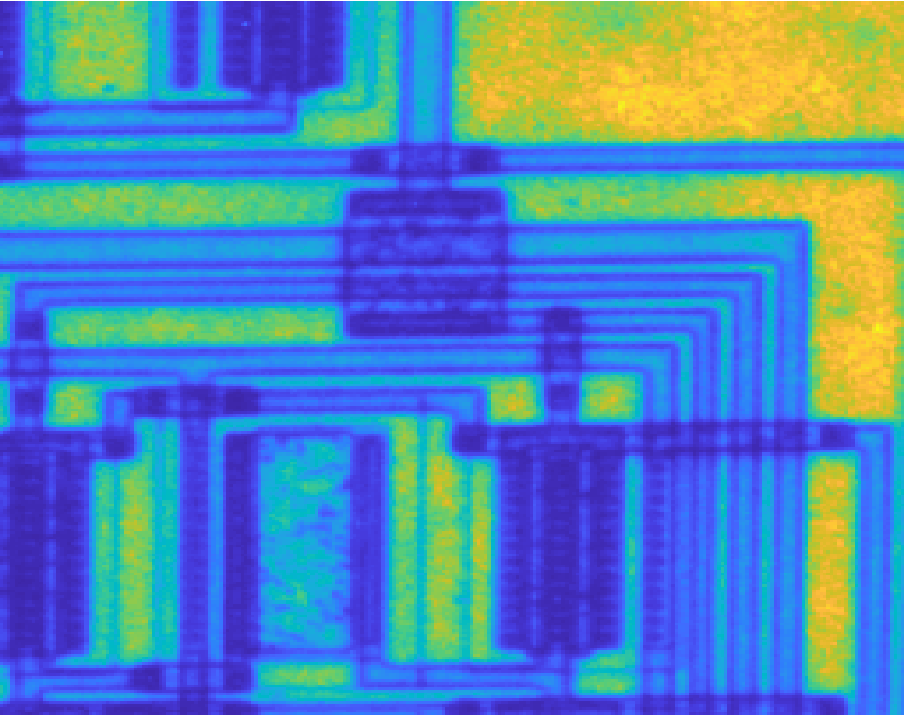} 
\end{tabular}
\begin{tabular}{ccc}
    {\scriptsize LSQR} &  {\scriptsize LSLU} & {\scriptsize sLSLU} \\ 
    \includegraphics[width=3.0cm]{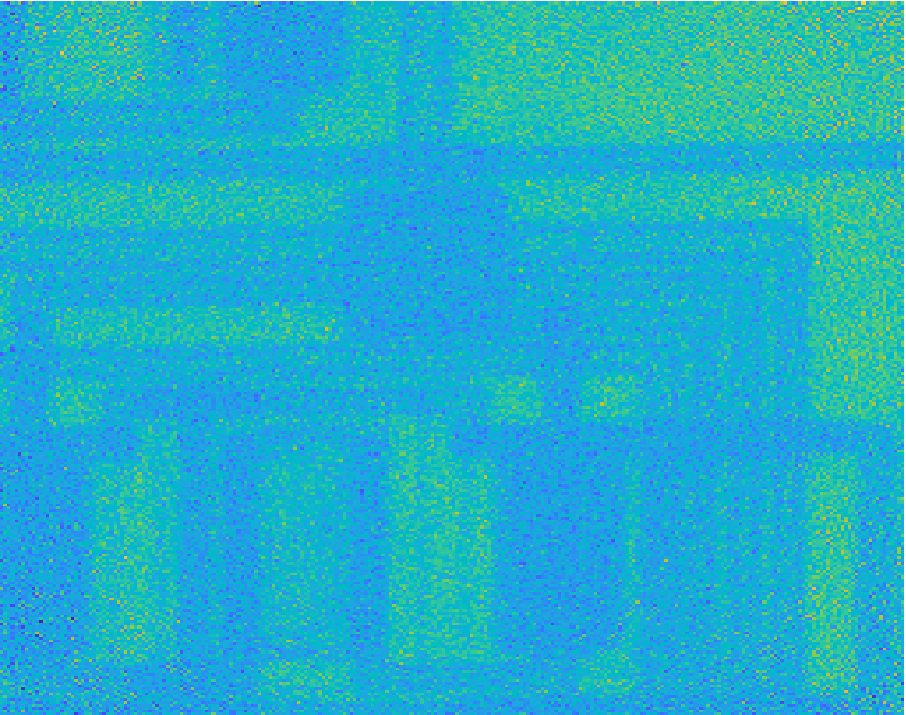} &  \includegraphics[width=3.0cm]{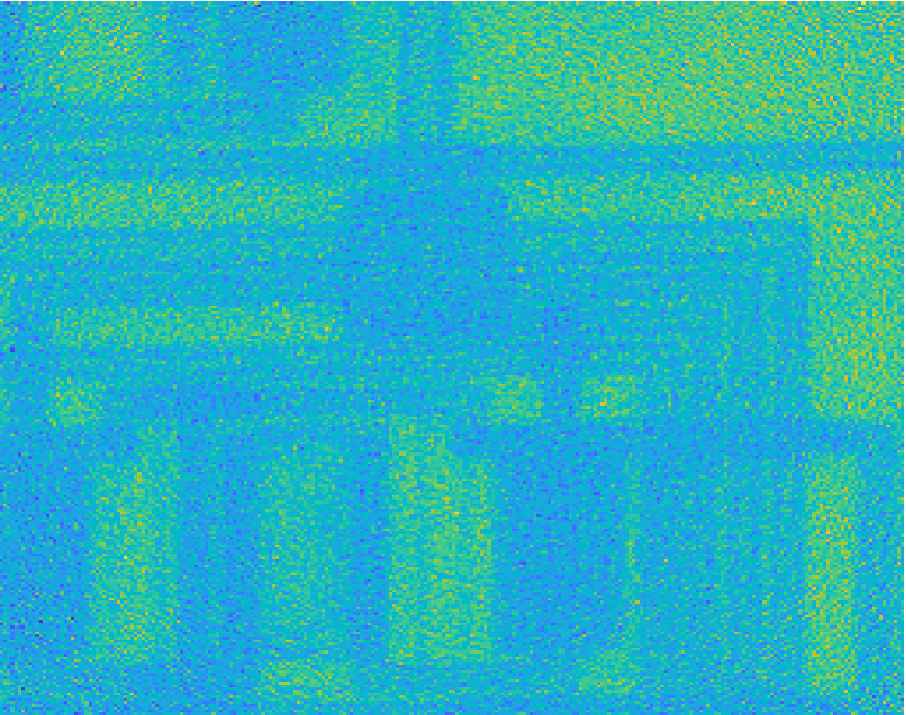} &
    \includegraphics[width=3.0cm]{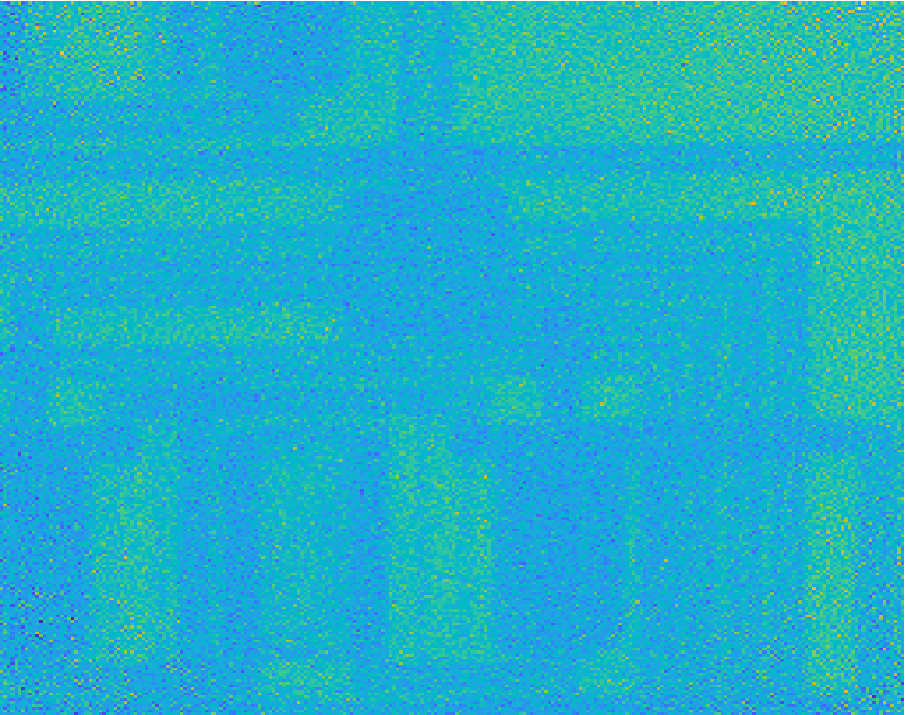}
\end{tabular} \\
\caption{Measured noisy data, true solution, and reconstructed images from LSQR, LSLU, and sLSLU. The image proportions are accurate but, to aid visualization, the relative size between images is not.}
\label{fig:rec_neutron_simulation}
\end{figure}

The residual norms are also plotted for comparison (see right plot in  \Cref{fig:sLSLU_rel}). In \Cref{fig:sLSLU_rel}, we find that residual norms for sLSLU closely follow the lower bound, which corresponds to residual norms for LSQR. Similar to the relative error plot, we find that the behavior of sLSLU aligns with LSQR, and the reconstructions are provided in \Cref{fig:rec_neutron_simulation}.

\begin{figure}[ht]
\centering
\begin{tabular}{cc}
    {\scriptsize Relative Reconstruction Error Norms } &  {\scriptsize Residual Norms} \\
     \includegraphics[width=0.45\textwidth]{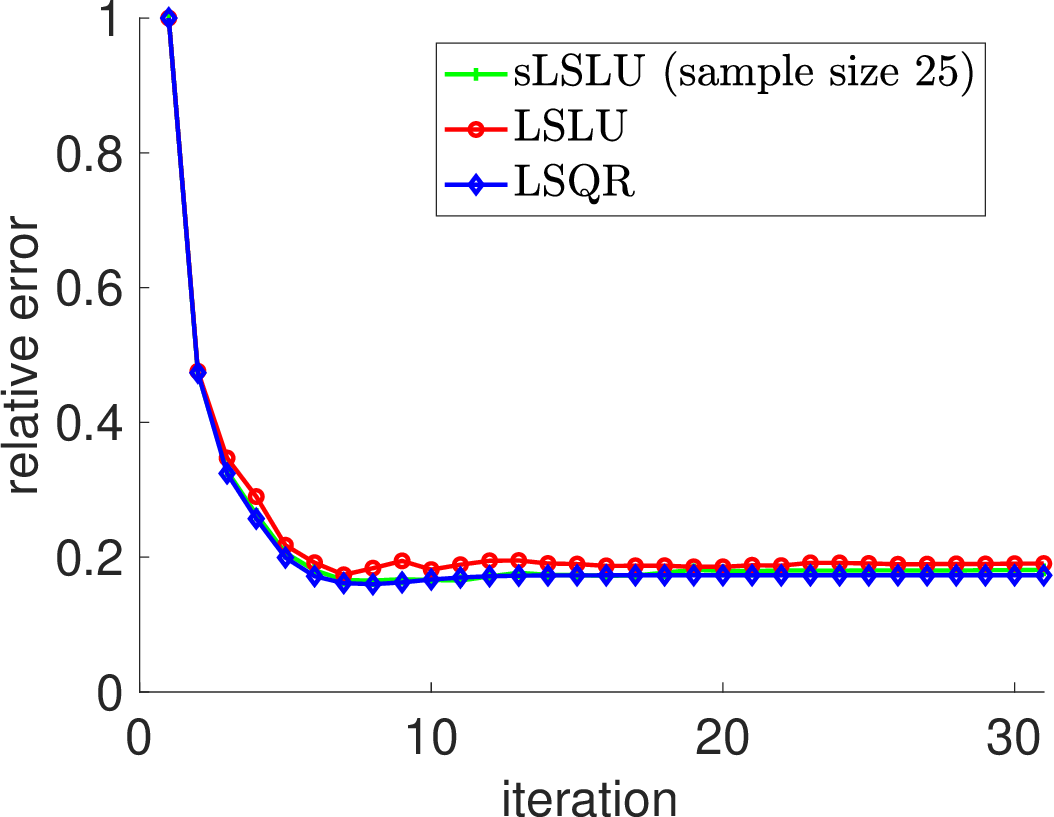} &  
   \includegraphics[width=0.45\textwidth] {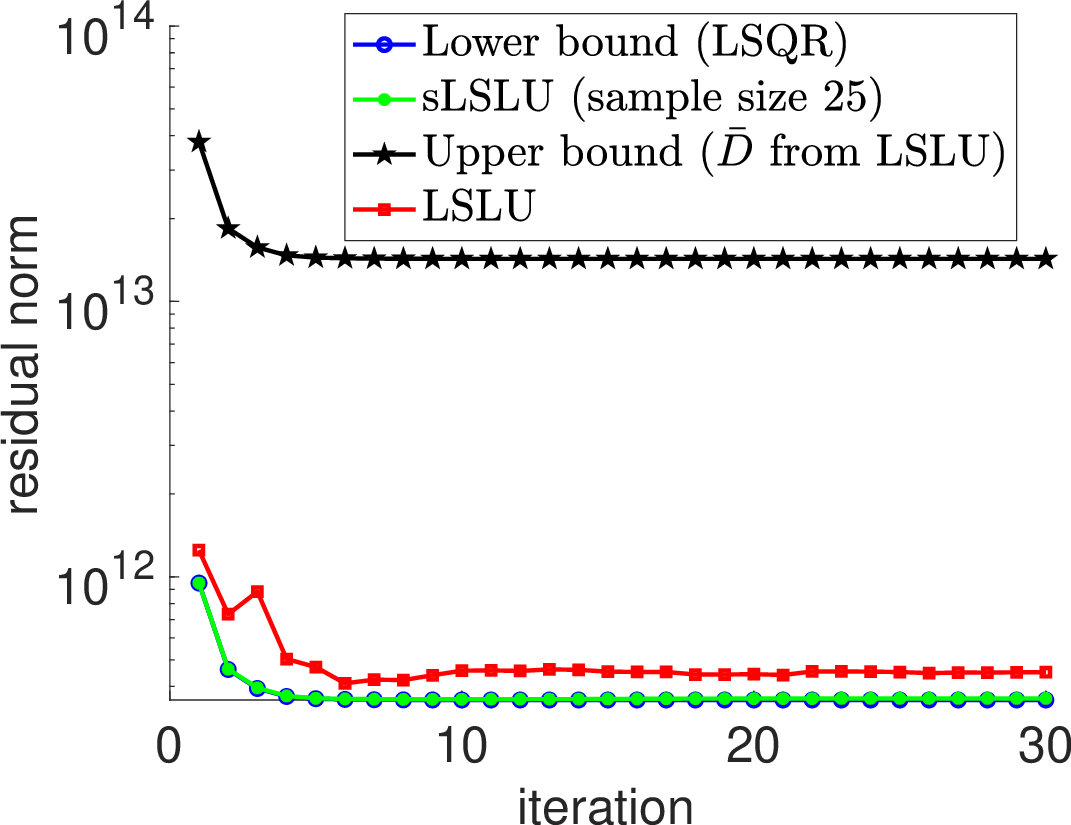} \\  
\end{tabular}
\caption{Relative reconstruction error norms (left) and residual norms (right). In this case, sketched LSLU with Tik. Reg. uses pivots dictated by the maximum absolute value from a set of randomly sampled coefficients (25).}
\label{fig:sHLSLU_rel}
\end{figure}

\begin{figure}[ht]
\centering
\title{Neutron Tomography Simulation} \\[12pt]
\begin{tabular}{cc}
    {\scriptsize Noisy Data} & {\scriptsize True Solution}\\
   
    \includegraphics[width=3.0cm]{noisy_neutron.eps} &
   \includegraphics[width=3.0cm]{true_solution_neutron.eps} \\
\end{tabular} 
\begin{tabular}{ccc}
{\scriptsize LSQR} &  {\scriptsize LSLU} & {\scriptsize sLSLU} \\ 
    \includegraphics[width=3.0cm]{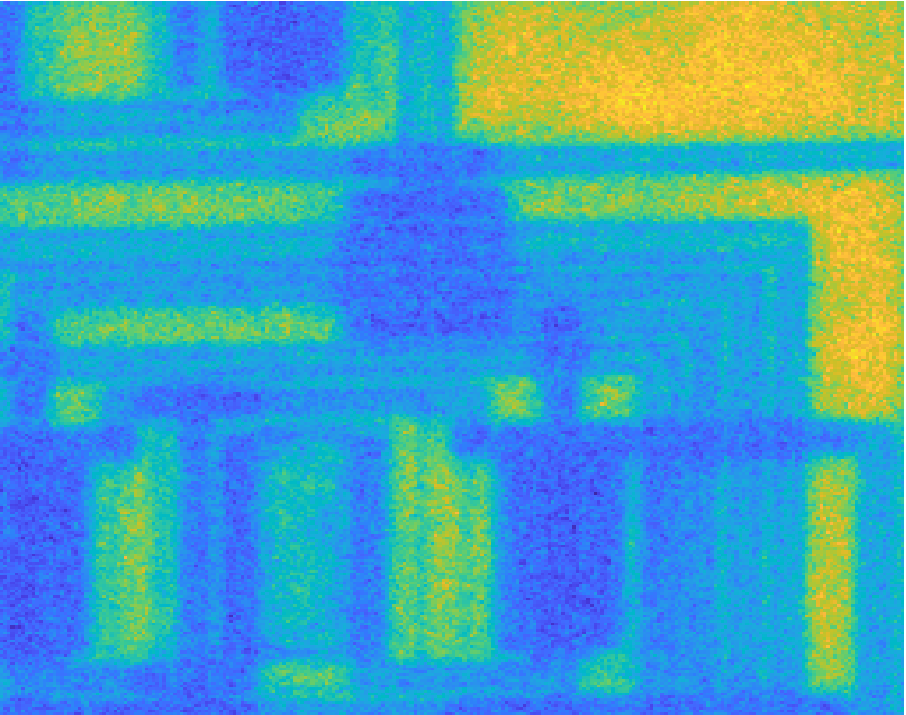} &  \includegraphics[width=3.0cm]{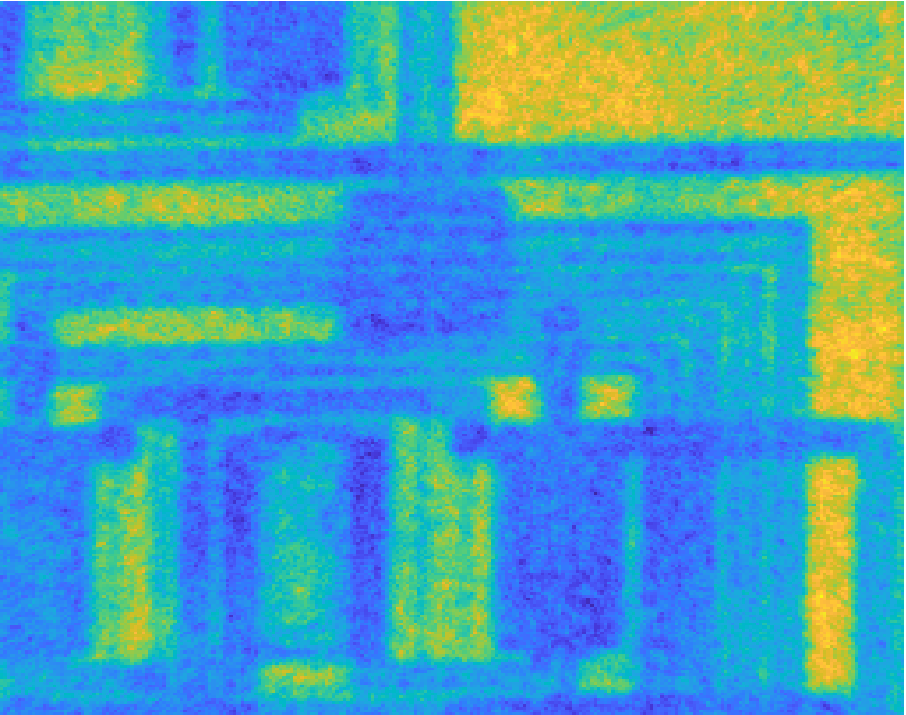} &
    \includegraphics[width=3.0cm]{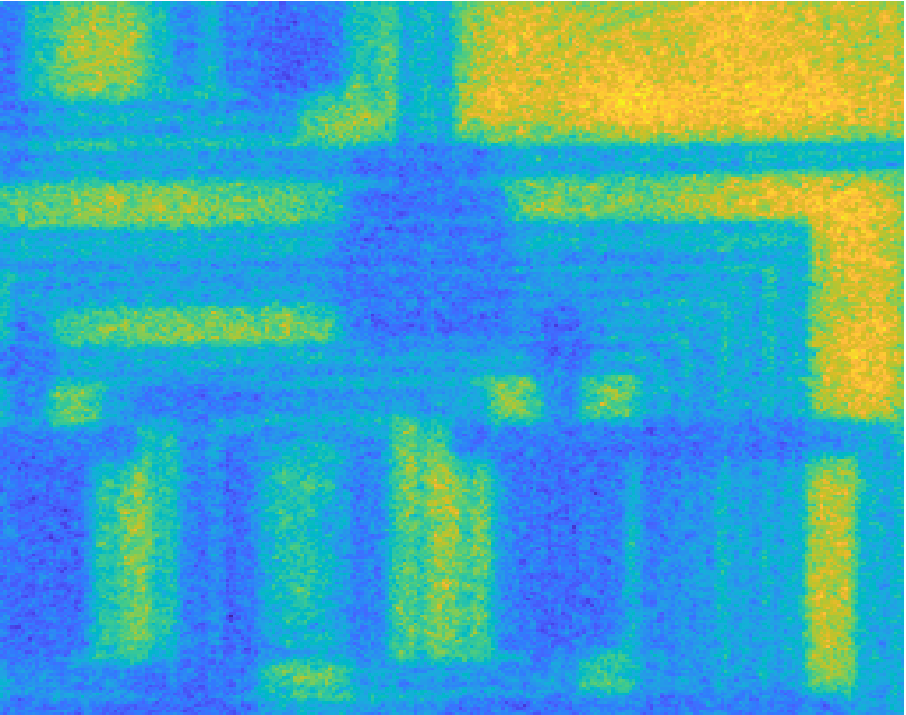} 
\end{tabular}\\
\caption{Measured noisy data, true solution, 
and reconstructed LSQR, LSLU, and sLSLU for solving the Tikhonov problem. The image proportions are accurate but, to aid visualization, the relative size between images is not.}
\label{fig:rec_neutron_simulation_hybrid}
\end{figure}

Next we consider the performance of sLSLU for the Tikhonov problem. We fix $\lambda = 26$, and we plot the residual norm of sLSLU compared to LSLU and LSQR in \Cref{fig:sHLSLU_rel}. We also provide the relative reconstruction error norms per iteration.  We observe that the inclusion of the regularization term stabilizes the semi-convergence for all methods.  For the Tikhonov problem, sLSLU mirrors the behavior of LSQR on the Tikhonov problem.  An adaptive approach to find a ``good" regularization parameter during the iterations is a topic of future work. Image reconstructions corresponding to $30$ iterations are provided in \Cref{fig:rec_neutron_simulation_hybrid}.

\subsection{Real data examples}
\begin{figure}[ht]
\centering
\title{Carved Cheese}
\begin{tabular}{cccc}
    {\scriptsize Noisy Data} &
    {\scriptsize LSQR} &  {\scriptsize LSLU} & {\scriptsize sLSLU}\\ 
    \includegraphics[width=2.9cm]{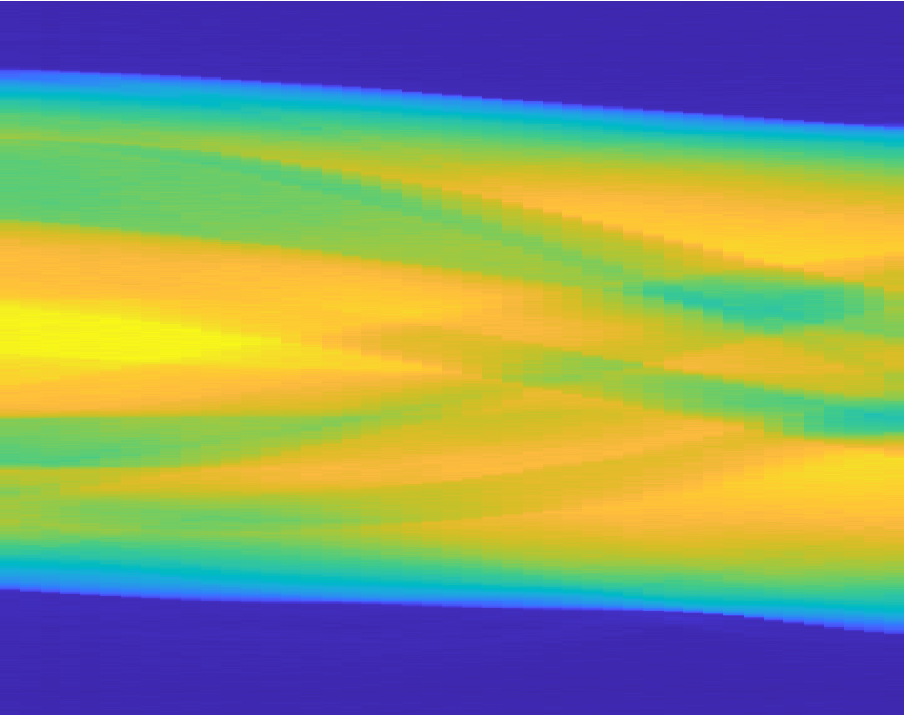} & 
    \includegraphics[width=2.9cm]{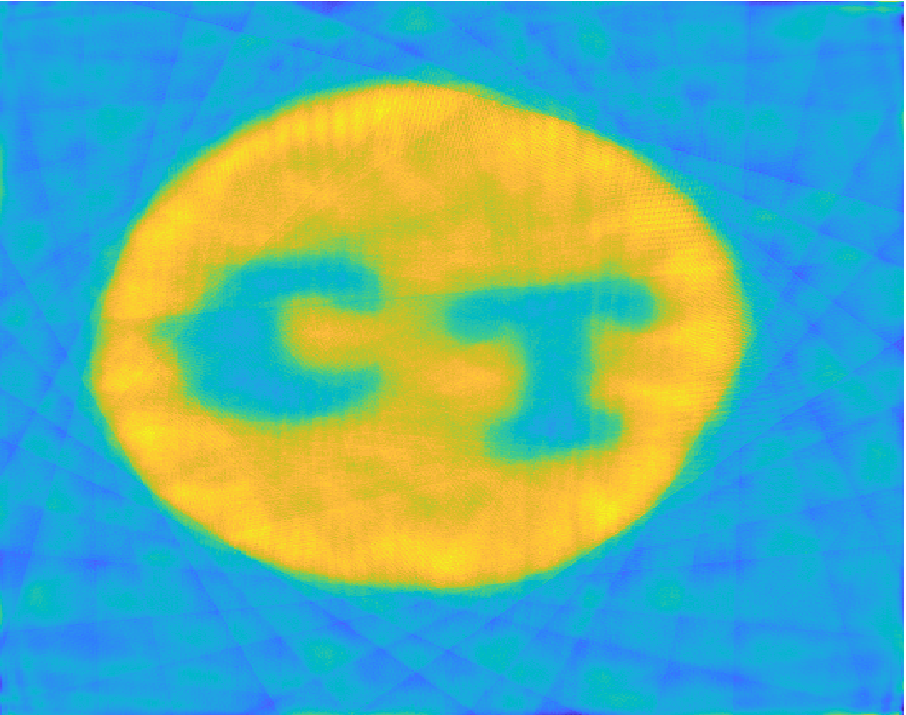} &  \includegraphics[width=2.9cm]{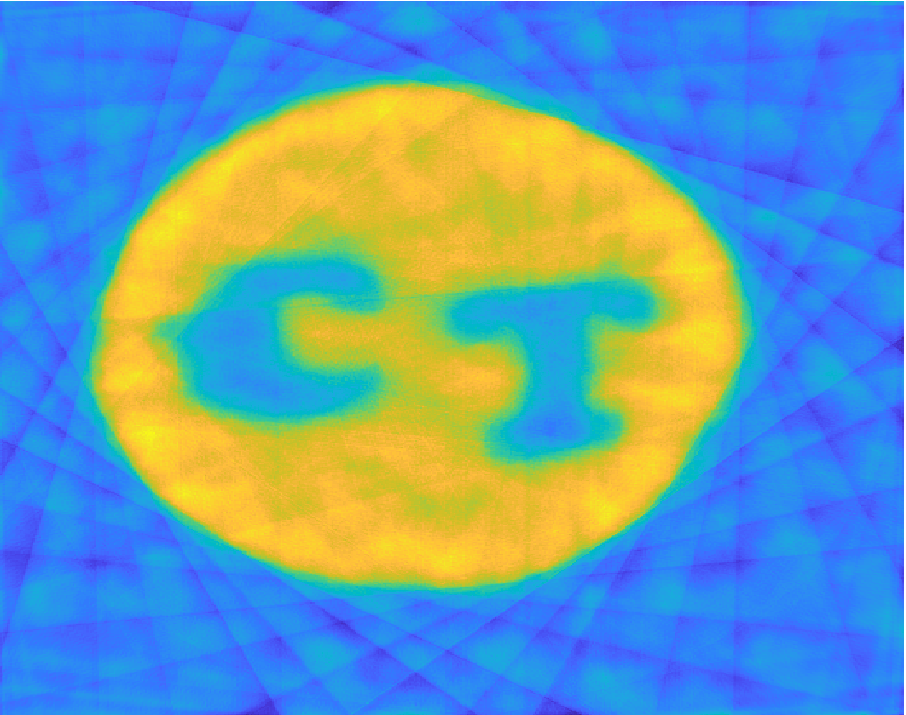} &
    \includegraphics[width=2.9cm]{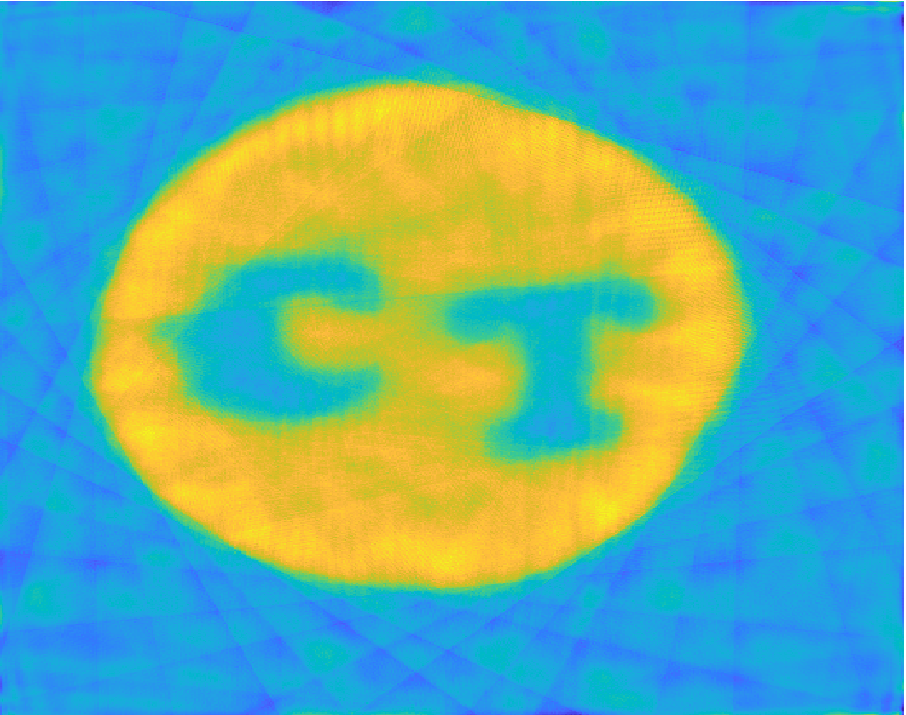} \\
\end{tabular} \\
\caption{Measured noisy data, and reconstructed images from LSQR, LSLU, and sLSLU. The image proportions are accurate but, to aid visualization, the relative size between images is not.}
\label{fig:rec_carved_cheese}
\end{figure}

\begin{figure}[ht]
\centering
\title{Walnut}\\
\begin{tabular}{cccc}
    {\scriptsize Noisy Data} &
    {\scriptsize LSQR} &  {\scriptsize LSLU} & {\scriptsize sLSLU}\\ 
    \includegraphics[width=2.9cm]{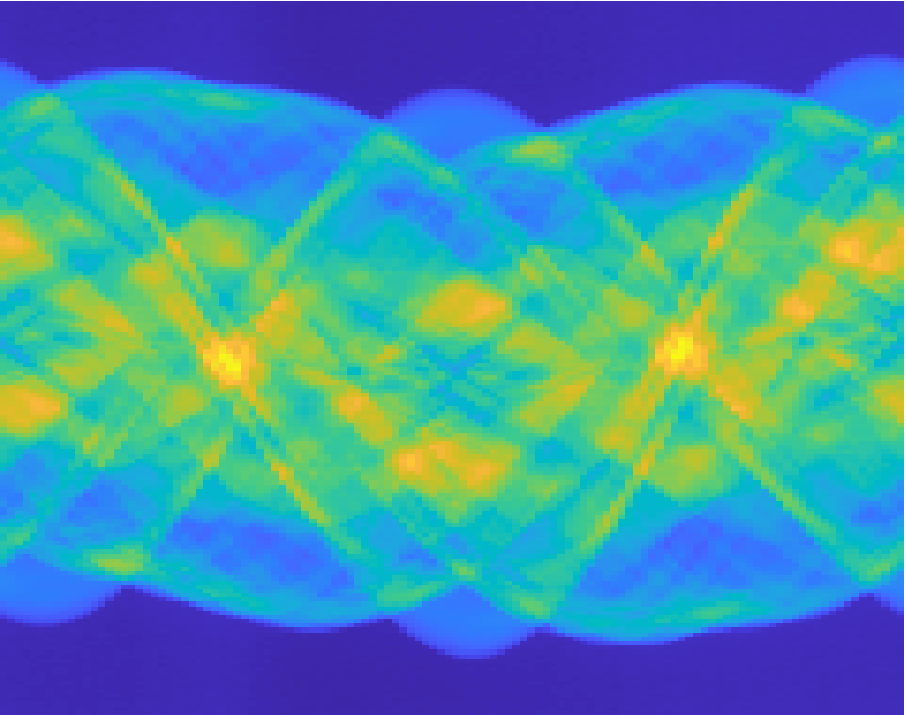} & 
    \includegraphics[width=2.9cm]{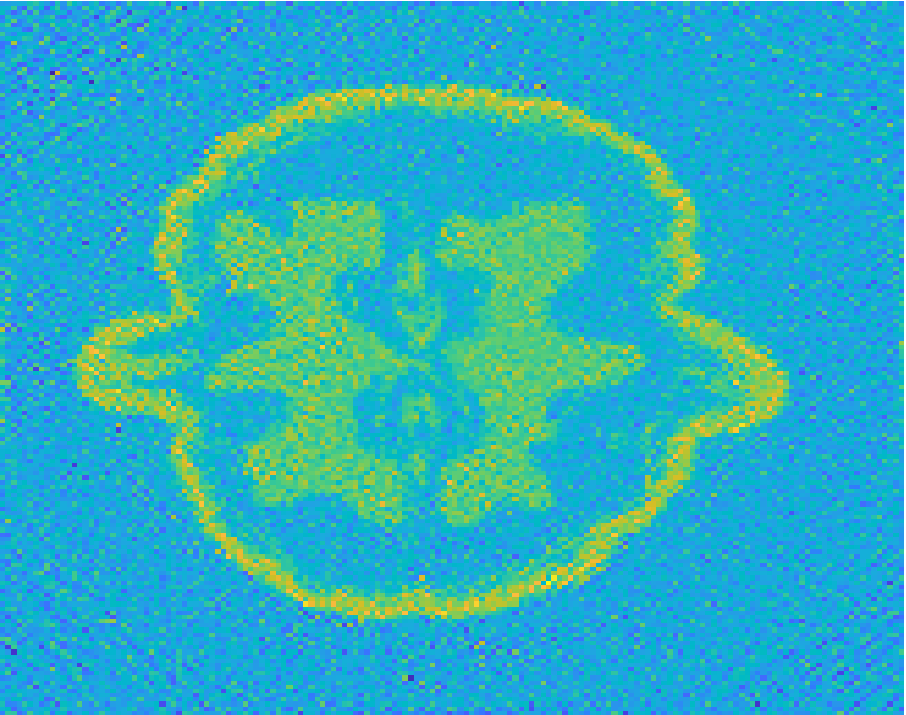} &  \includegraphics[width=2.9cm]{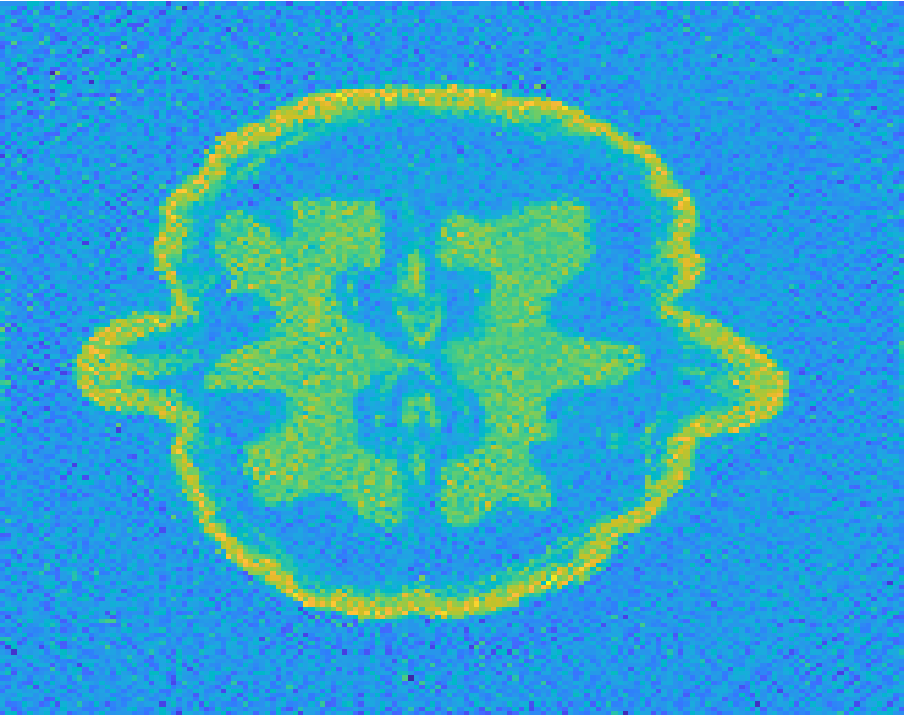} &
    \includegraphics[width=2.9cm]{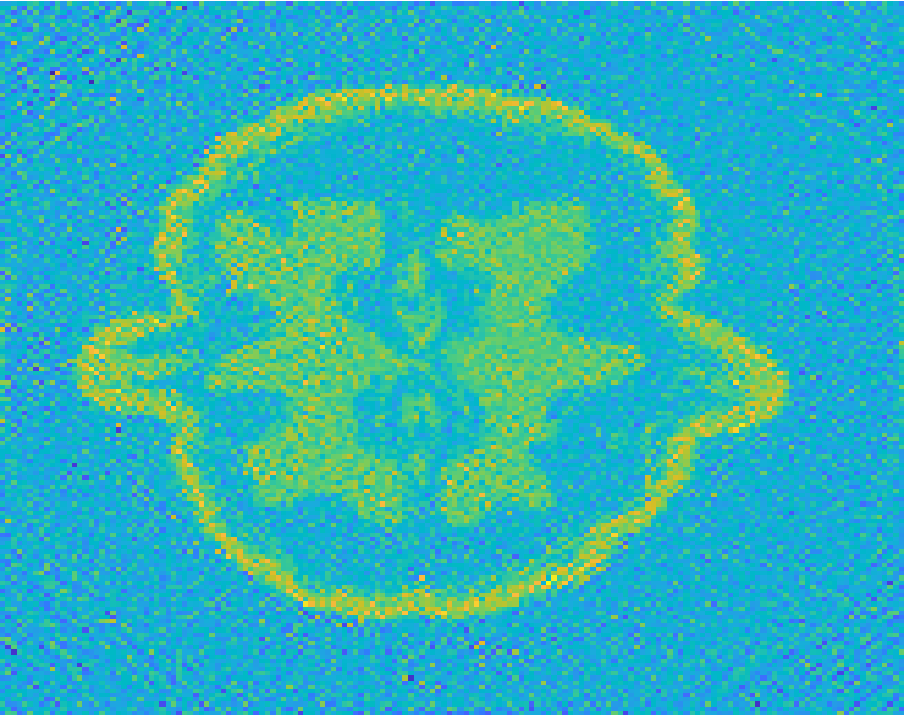} \\
\end{tabular} \\

\caption{Measured noisy data, and reconstructed images from LSQR, LSLU, and sLSLU. The image proportions are accurate but, to aid visualization, the relative size between images is not.}
\label{fig:rec_walnut}
\end{figure}

The Finnish Inverse Problems Society has provided the following open access datasets: a tomographic x-ray of carved cheese and a walnut. Both datasets consist of X-ray sinograms where each sinogram is obtained by fan-beam projection. The observed data for the carved cheese dataset containing $360$ projections and the walnut dataset containing $120$ projections are provided in Figures \ref{fig:rec_carved_cheese} and \ref{fig:rec_walnut} respectively. For these problems, there is no true image, so we rely on residual norms per iteration to compare algorithms. 

\begin{figure}[ht]
\centering
\begin{tabular}{cc}
    {\scriptsize Carved Cheese} &  {\scriptsize Walnut} \\ 
     \includegraphics[width=0.45\textwidth]{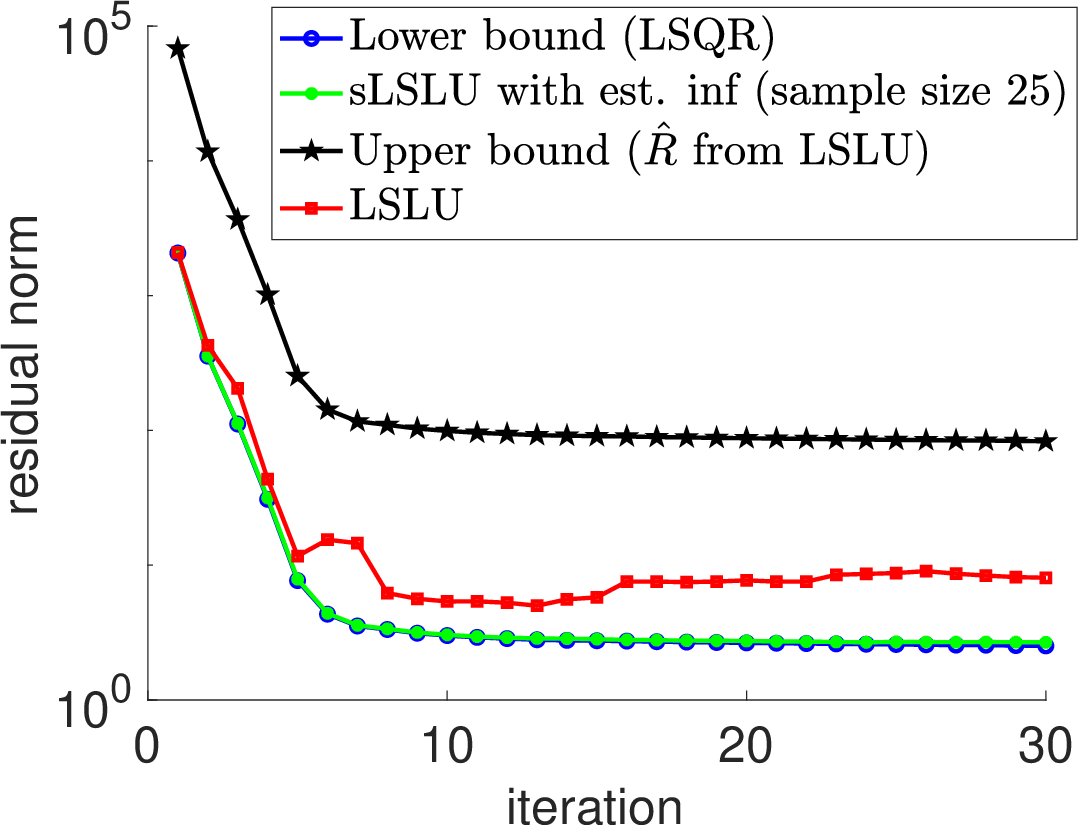}  &
   \includegraphics[width=0.45\textwidth]{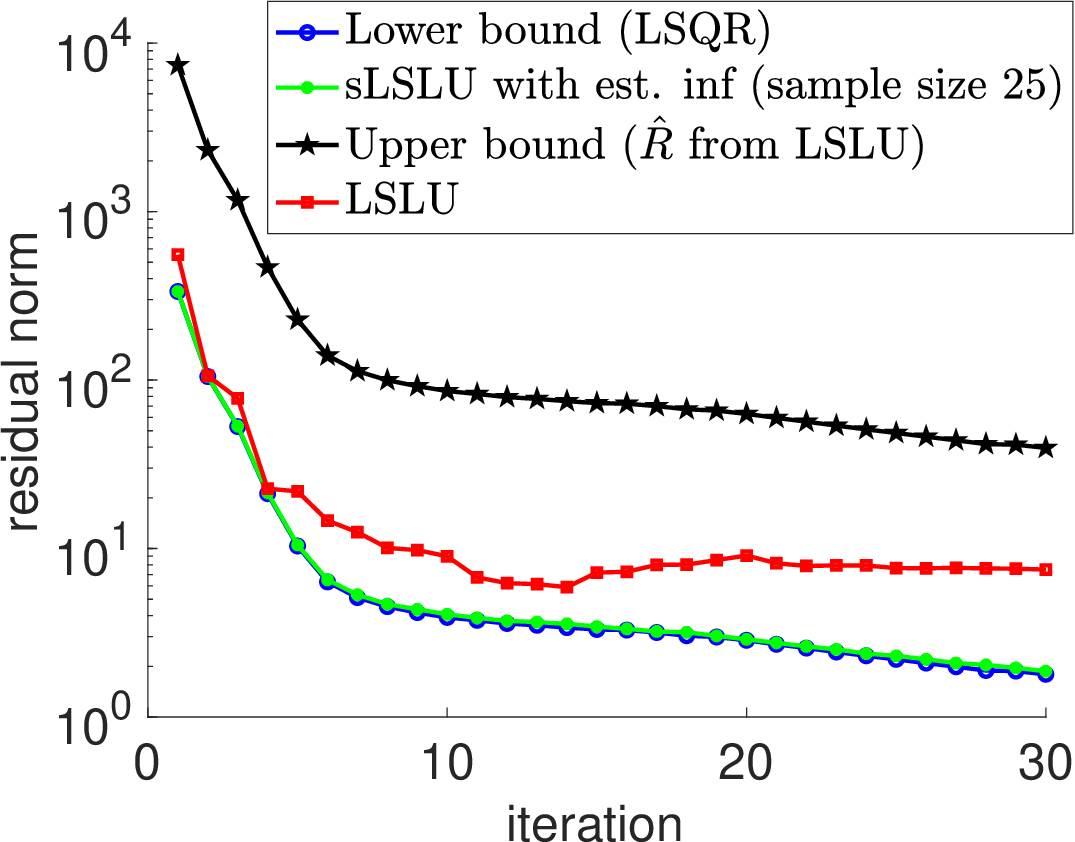} \\  
    \end{tabular}   
\caption{Residual norms per iteration for sLSLU and LSLU, as well as corresponding bounds from Theorem 2.1 of \cite{brown2024hlslu}. Note that the lower bound corresponds to LSQR residual norms.}
\label{fig:sLSLU residuals}
\end{figure}

To illustrate the behavior of the residual norms for sLSLU, LSLU, and LSQR as well as the bounds in Theorem $2.1$ of \cite{brown2024hlslu}, we plot in \Cref{fig:sLSLU residuals} the residual norms per iteration for the carved cheese and walnut datasets. The samples size to approximate the infinity norm contains $25$ entries. For both datasets, we observe that the residual norms for sLSLU and LSQR remain close together, with the residual norms for LSLU being a bit larger. Therefore, we may expect that the approximate solutions from sLSLU should mirror those from LSQR.  This is verified in \Cref{fig:rec_carved_cheese} and \Cref{fig:rec_walnut}, where the reconstructed images using LSQR, LSLU, and sLSLU are provided. All reconstructions correspond to iteration $30$.

\begin{figure}[ht]
\centering
\begin{tabular}{cc}
    {\scriptsize Carved Cheese} &  {\scriptsize Walnut} \\ 
      \includegraphics[width=0.45\textwidth]{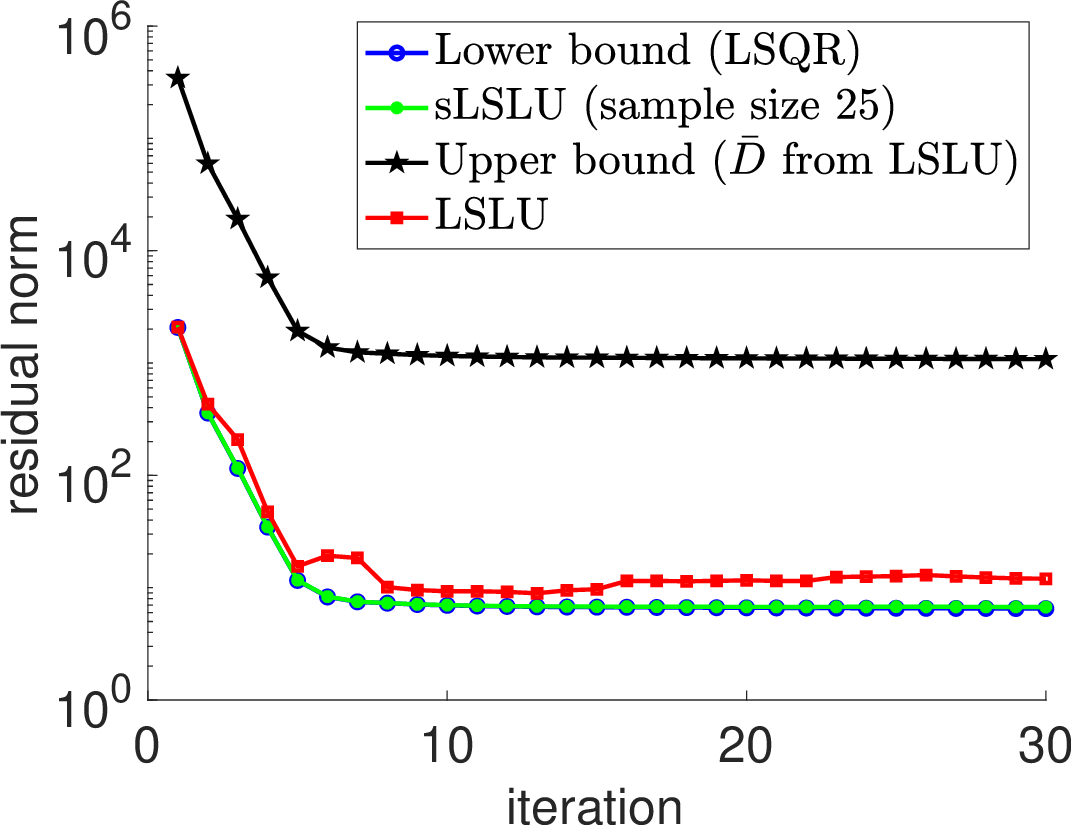} &
     \includegraphics[width=0.45\textwidth]{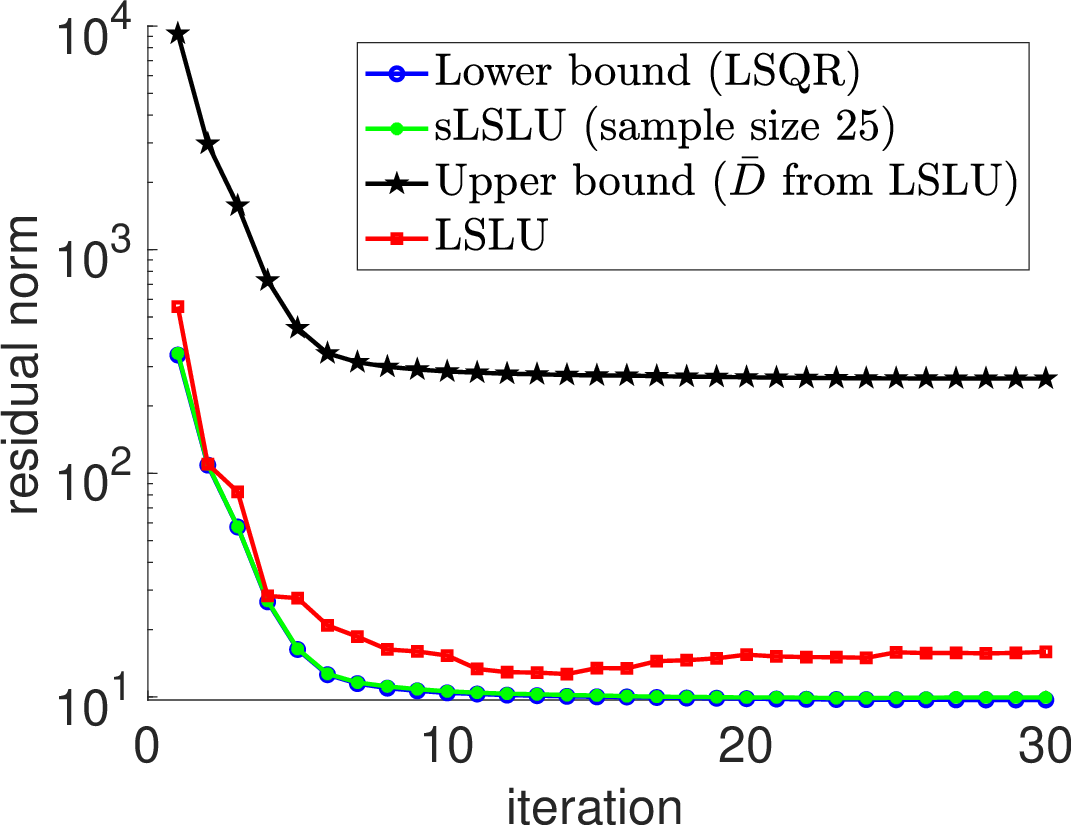} \\
    \end{tabular}
\caption{Residual norms per iteration for the Tikhonov problem correspond to sLSLU, LSLU, as well as corresponding bounds from Theorem 3.1 of \cite{brown2024hlslu}. Note that the lower bound corresponds to LSQR. The regularization parameter $\lambda=1$.}
\label{fig:hybrid_residual}
\end{figure}

\begin{figure}[ht]
\centering
\title{Carved Cheese}\\
\begin{tabular}{cccc}
    {\scriptsize Noisy Data} &
    {\scriptsize  LSQR} &  {\scriptsize  LSLU} & {\scriptsize sLSLU }\\ 
    \includegraphics[width=2.9cm]{carved_cheese_noisy.eps} & 
    \includegraphics[width=2.9cm]{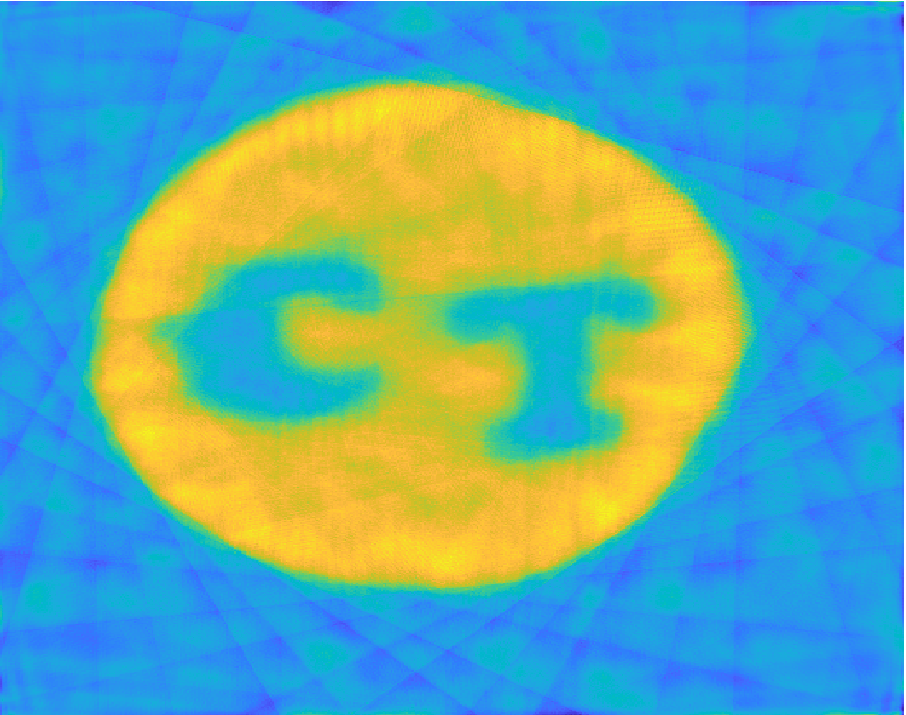} &  \includegraphics[width=2.9cm]{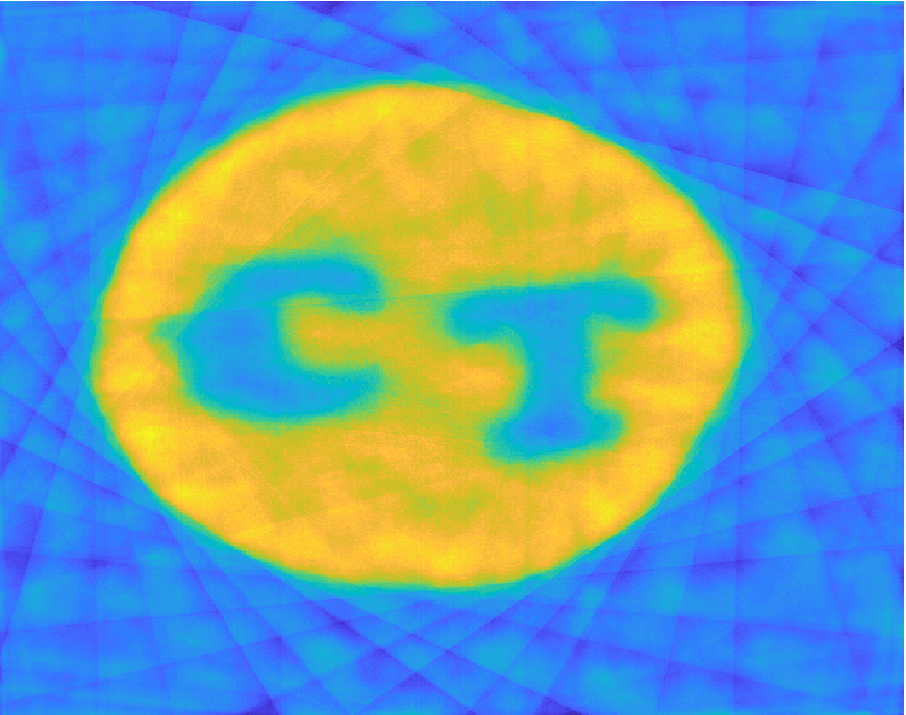} &
    \includegraphics[width=2.9cm]{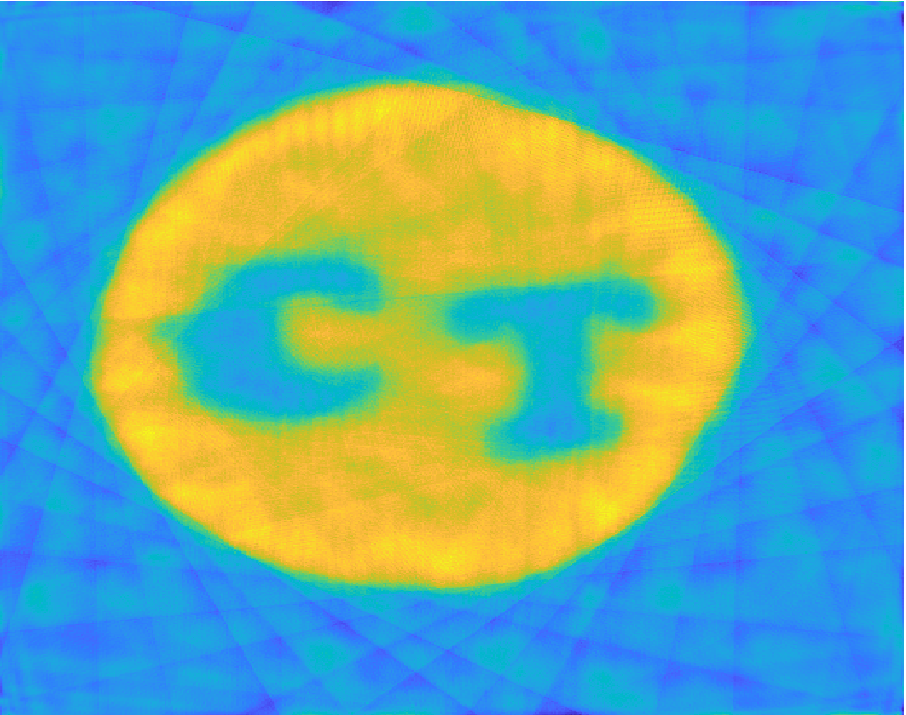} \\
\end{tabular} \\

\caption{Measured noisy data, and Tikhonov reconstructions for LSQR,  LSLU, and sLSLU. The image proportions are accurate but, to aid visualization, the relative size between images is not.}
\label{fig:rec_carved_cheese_with_tik}
\end{figure}

\begin{figure}[ht]
\centering
\title{Walnut}\\
\begin{tabular}{cccc}
    {\scriptsize Noisy Data} &
    {\scriptsize  LSQR} &  {\scriptsize LSLU} & {\scriptsize sLSLU }\\ 
    \includegraphics[width=2.9cm]{walnut_noisy.eps} & 
    \includegraphics[width=2.9cm]{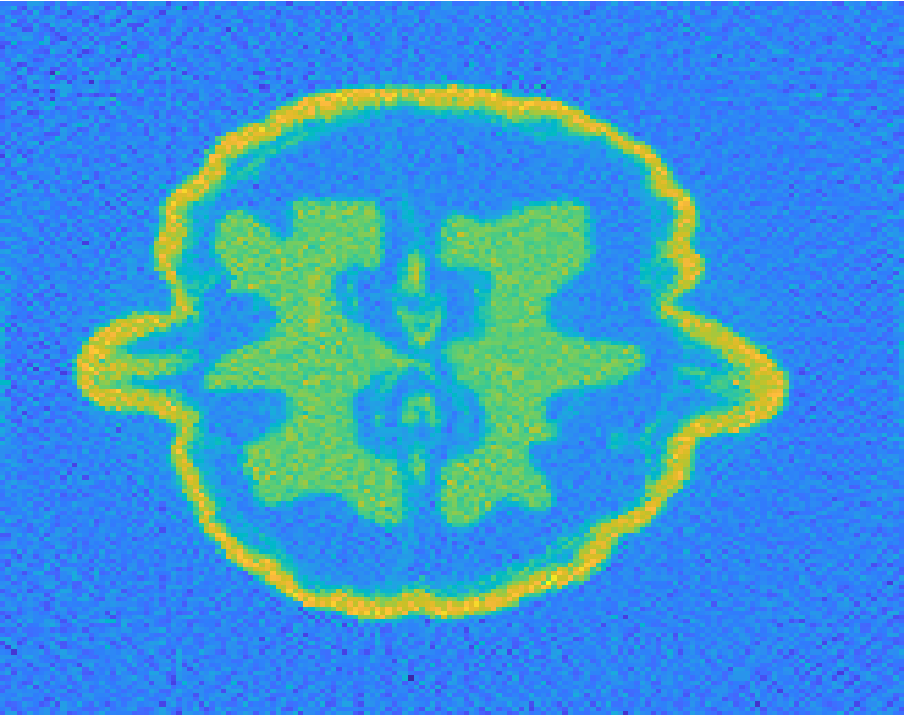} &  \includegraphics[width=2.9cm]{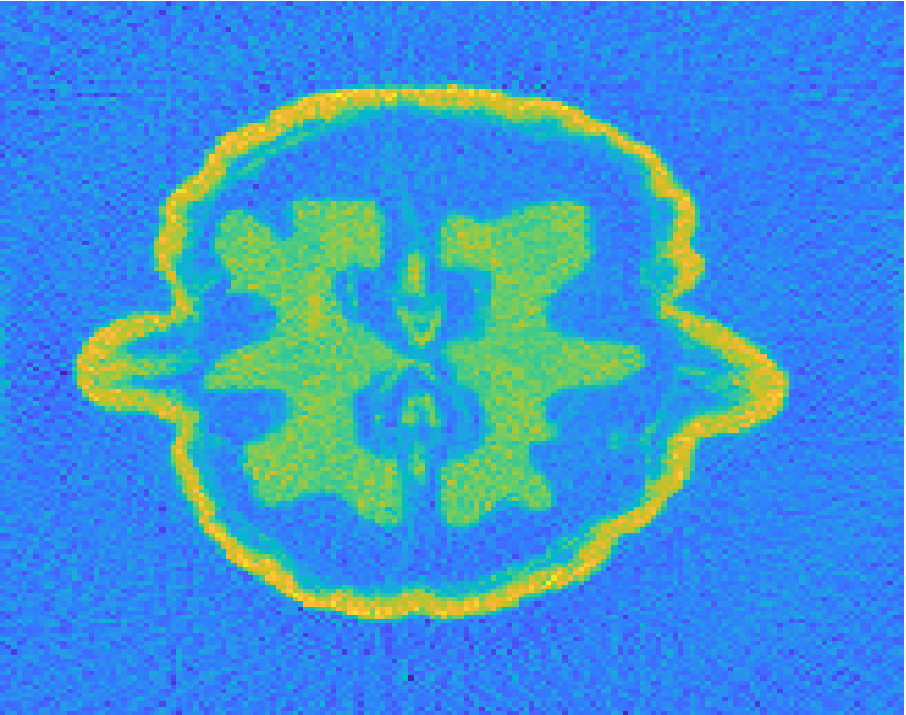} &
    \includegraphics[width=2.9cm]{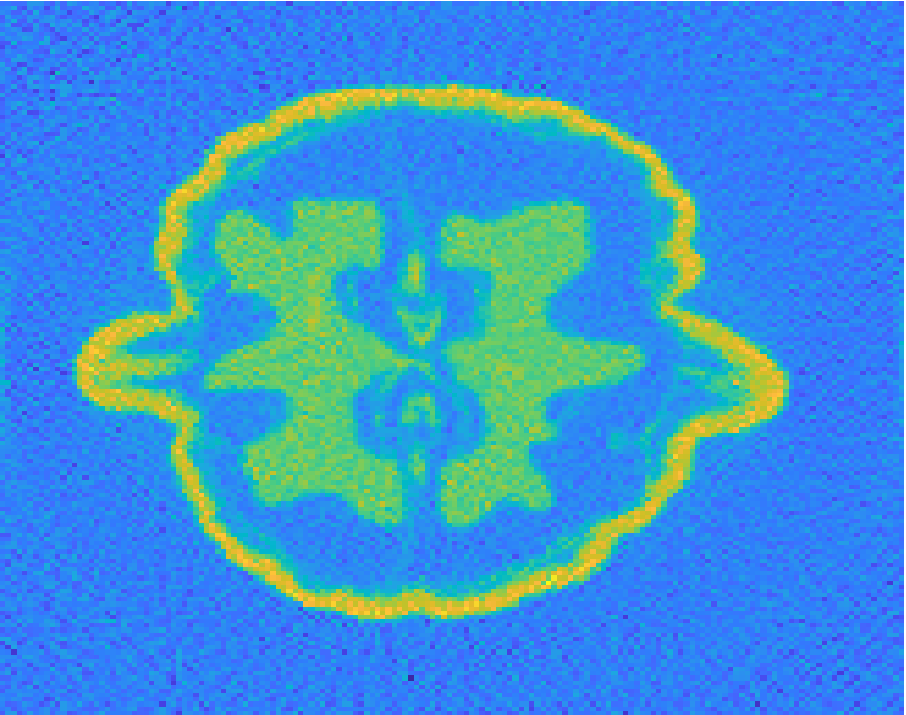} \\
\end{tabular} \\

\caption{Measured noisy data, and Tikhonov reconstructions for LSQR,  LSLU, and sLSLU. The image proportions are accurate but, to aid visualization, the relative size between images is not.}
\label{fig:rec_walnut_with_tik}
\end{figure}

 Finally, we consider sLSLU with Tikhonov regularization for these examples, where we plot the the residual norm of sLSLU with Tikhonov regularization in \Cref{fig:hybrid_residual}. We fix $\lambda = 1$ for both problems. Similar to the nonregularized problems, the residual norms for sLSLU closely follow the lower bound, which corresponds to residual norms for LSQR. Thus, we may expect that provided we have a ``good" estimate for the regularization parameter, sLSLU with Tikhonov regularization will produce a better approximation of the solution than LSLU on the Tikhonov problem. We also provide reconstructed images for both datasets in \Cref{fig:rec_carved_cheese_with_tik} and \Cref{fig:rec_walnut_with_tik}.

\backmatter

\section{Conclusions} \label{sec:conc}
In this paper, we introduce two new inner-product free Krylov methods, sCMRH and sLSLU, that incorporate randomization techniques for solving large-scale linear inverse problems. Both methods are based on the Hessenberg method with partial pivoting for building bases that span Krylov subspaces, and hence do not require inner-product computations (e.g., orthogonalizations). Also, both sCMRH and sLSLU exploit randomized sketching to solving the projected problems, thereby producing solutions with a smaller residual norm compared to existing inner-product free Krylov methods. Numerical experiments show that the performance of sCMRH is comparable to that of GMRES and the perforamcne of sLSLU is comparable to that of LSQR.  Moreover, sCMRH and sLSLU have smaller residual norm solutions, compared to CMRH and LSLU respectively. The sketched Krylov methods can be adapted to incorporate Tikhonov regularization provided that an appropriate regularization parameter is selected. Since sCMRH and sLSLU are all inner-product free, they may be useful in solving problems with mixed-precision and parallel computing, which is a topic of future work.

\bmhead{Acknowledgements}
We would like to acknowledge Ethan Epperly's blog, which provided us with nice expositions:
\href{https://www.ethanepperly.com/index.php/blog}{https://www.ethanepperly.com/index.php/blog}\\

\bmhead{Funding}
This work was partially funded by the U.S. National Science Foundation, under grants DMS-2038118, DMS-2411197 and DMS-2208294. Any opinions, finding, and conclusions or recommendations expressed in this material are those of the author(s) and do not necessarily reflect the views of the National Science Foundation. 

\bibliography{references}%

\end{document}